\documentclass[10pt]{amsart}
\usepackage{hyperref}
\usepackage{amssymb, amsmath, amsthm}
\newtheorem{thm}{Theorem}[section]
\newtheorem{lemma}{Lemma}[section]
\newtheorem{cor}{Corollary}[section]
\newtheorem{prop}{Proposition}[section]
\newtheorem{defn}{Definition}[section]

\theoremstyle{definition}
\newtheorem{rem}{Remark}[section]

\newcommand{\PP}{\mathcal{P}}  
\newcommand{\w}{\omega}        
\newcommand{\ed}{\mathbf{d}}   
\newcommand{\J}{\mathbf{J}}    
\newcommand{\D}{\mathbf{D}}    

\newcommand{\R}{\mathbb{R}}    
\newcommand{\Ad}{\mathrm{Ad}}  
\newcommand{\ad}{\mathrm{ad}}  
\newcommand{\FL}{\mathbb{F}\mathrm{L}}  
\newcommand{\Sl}{\mathbf{S}}   
\newcommand{\OO}{\mathcal{O}}  
\newcommand{\g}{\mathfrak{g}}
\newcommand{\h}{\mathfrak{h}}
\newcommand{\lir}{\mathfrak{r}}
\newcommand{\q}{\mathfrak{q}}
\newcommand{\be}{\begin{equation}}
\newcommand{\ee}{\end{equation}}
\newcommand{\bea}{\begin{eqnarray}}
\newcommand{\eea}{\end{eqnarray}}
\newcommand{\II}{\mathbb{I}} 
\newcommand{\rII}{\hat{\mathbb{I}}_0} 
\newcommand{\Proj}{\mathbb{P}} 
\newcommand{\corr}{\mathrm{corr}\,} 
\newcommand{\Hor}{\mathrm{Hor}} 
\newcommand{\Ar}{\mathrm{Ar}} 
\newcommand{\restr}[1]{\vrule height3ex width.4pt depth1.4ex\lower1.4ex\hbox{\scriptsize $\,#1$}}
\newcommand{\rrestr}[1]{\vrule height2ex width.4pt depth0.9ex\lower0.9ex\hbox{\scriptsize $\,#1$}}

\makeatletter \@addtoreset{figure}{section}
\def\thefigure{\thesection.\@arabic\c@figure}
\def\fps@figure{h, t}
\@addtoreset{table}{bsection}
\def\thetable{\thesection.\@arabic\c@table}
\def\fps@table{h, t}
\@addtoreset{equation}{section}

\makeatother

\allowdisplaybreaks
\def\intprod{\mathbin{\hbox to 6pt{%
                    \vrule height0.4pt width5pt depth0pt
                    \kern-.4pt
                    \vrule height6pt width0.4pt depth0pt\hss}}}

\begin{document}

\title[Stability in Simple Mechanical Systems]{Stability of Relative
Equilibria with Singular Momentum Values in Simple Mechanical Systems}

\author{Miguel Rodr\'{\i}guez-Olmos}
\thanks{corresponding email: miguel.rodriguez@epfl.ch,
phone: +41 21 6935507, fax: +41 21 6935839}
\address{Ecole Polytechnique F\'ed\'erale de Lausanne (EPFL), Section de
Math\'ematiques. CH-1015 Lausanne, Switzerland.}

\begin{abstract}A method for testing stability of
relative equilibria in Hamiltonian systems of the form ``kinetic +
potential energy" is presented. This method extends previously
existing techniques to the case of non-free group actions and
singular momentum values. A normal form for the symplectic
matrix at a relative equilibrium is also obtained.
\end{abstract}

\maketitle

\section{Introduction}
A relative equilibrium in a dynamical system with symmetry is a
point in phase space for which its dynamical evolution is
contained in a group orbit. The study of relative equilibria in
symmetric Hamiltonian systems has been around for a long time,
with its origins in the field of analytical dynamics  and more
recently using the modern symplectic and Poisson geometric setup.
Relative equilibria are important since they are the analogues to
equilibrium states in systems with symmetry, formalized with the
action of a Lie group on the phase space. In physical
applications, the only observable equilibrium states are those
which are stable under small perturbations. Similarly, in the
symmetric context, the only observable relative equilibria are
those which are stable in some adequate sense. Based in Noether's
Theorem, geometrized in the property of the invariance of the
level sets of the momentum map, the notion of stability generally
adopted is that of \emph{$G_\mu$-stability} introduced in
\cite{Patrick}, and that is closely related to the Lyapunov
stability of the induced Hamiltonian flow on the reduced phase
space.

In the field of analytical dynamics, the classical Routh Theorem
gives conditions on the stability of steady motions which keep
stationary the value of a first integral of a dynamical system for
fixed values of the others (see for instance \cite{Routh} and the
treatments based on the Routh Theorem in \cite{Kara, Pascal}).
Relative equilibria are seen in this context as steady motions for
systems having cyclic coordinates due to the existence of a
symmetry group, for which the components of the momentum map
together with the energy provide a set of first integrals. In the
last decades, the implementation of these principles within the
field of Geometric Mechanics has been studied. This has produced
methods (like the Energy-Momentum Method \cite{Patrick} and the
Energy-Casimir Method \cite{Arnold2}, see also \cite{MarLec} for
an overview) to test the stability of relative equilibria in
Hamiltonian systems for arbitrary symmetry groups and momentum
values. These methods exploit Noether's Theorem and the symplectic
and Poisson geometry of the phase space. In the case that the
relative equilibrium under study lies in a regular value $\mu$ of
the momentum map, the Energy-Momentum Method of \cite{Patrick}
provides a technique to test its stability modulo the action of
$G_\mu$, the stabilizer of the momentum value $\mu$ under the
coadjoint representation of $G$. This was generalized in
\cite{LeSi} and \cite{OrRa99} to also cover the case when the
momentum value is singular (this happens when the symmetry group
does not act freely on phase space), assuming $G_\mu$  is compact.
Also, in \cite{Montaldi} and \cite{PatRobWul} stability of
relative equilibria satisfying several other hypotheses is
investigated.

A very important particular kind of Hamiltonian system is the
class of \emph{simple mechanical systems}, paradigmatic of
Classical Mechanics, since many Hamiltonian system of physical
interest lie in this category or can be obtained from a simple
mechanical system by some suitable reduction process. These have
as phase space the cotangent bundle of a Riemannian manifold
(called configuration space) equipped with its canonical
symplectic form, and the Hamiltonian function is of the type
``kinetic plus potential'' energy, where the kinetic energy is
given by the norm obtained from the Riemannian metric, and the
potential energy is the pullback of a function defined on the
base. Symmetry in this systems is implemented by the lift of an
isometric action on the base that preserves the potential energy.
This big amount of extra structure with respect to general
Hamiltonian systems on arbitrary symplectic or Poisson phase
spaces implies that in simple mechanical systems everything is
constructible from the knowledge of the configuration space, its
Riemannian structure, the action of the symmetry group on it and
the choice of an invariant potential energy. Therefore, it is
reasonably to expect that the stability methods referred
previously will particularize in a way that the involved
computational complexity will simplify considerably, in that one
would work at the level of the configuration space, instead of on
its twice dimensional cotangent bundle.

In the case of regular relative equilibria, this refinement of the
Energy-Momentum Method has been worked out in \cite{SiLeMar}, and
the obtained stability test for relative equilibria in simple
mechanical systems is known as reduced Energy-Momentum Method. Its
conditions for $G_\mu$-stability are reduced from the level of
phase space to the level of configuration space. This method has
the highest degree of sophistication among the different stability
tests available in the literature of symmetric Hamiltonian
systems, and as part of it, it provides a block-diagonalization
technique that allows to express the linearization of the
Hamiltonian vector field at a relative equilibrium in a way
adapted both to the symmetry of the system and to the fibered
structure of the phase space. This block-diagonalization yields
also further simplifications in the stability analysis.

Surprisingly, in the very frequent and important case of singular
momentum values such a refinement for simple mechanical systems
has not been studied in detail, and thus the application range of
the reduced Energy-Momentum Method is severely limited. Indeed,
the literature of applications of the theory of relative
equilibria is full of examples in which singular relative
equilibria of simple mechanical systems are studied with general
geometric and Hamiltonian techniques which neglect their extra
structure, in particular for the stability analysis. This paper
provides a solution to this situation by obtaining a
generalization to singular momentum values of the reduced
Energy-Momentum Method and its main features.

In Section 2 we quickly review the theory of relative equilibria
for general Hamiltonian systems and simple mechanical systems, and
we collect some of the standard results on their stability by
geometric methods. Section 3 is a necessary technical interlude on
the properties and geometry of a distinguished symplectic
component of the linear slice for a cotangent-lifted action, and
most of our subsequent results will rely on this section. In
Section 4 our main result, Theorem \ref{REM}, is stated, providing
an extension of the reduced Energy-Momentum Method of
\cite{SiLeMar} applicable to relative equilibria with singular
momentum values. Section 5 applies this result to a classical
example of a relative equilibrium with a singular momentum value
in a well-known simple mechanical system consisting of an
axisymmetric rigid body with a fixed point in an homogeneous
gravity field. It is shown how the application of our method
simplifies the stability analysis with respect to the application
of the methods developed for general Hamiltonian systems. In
Section 6 we extend the block-diagonalization result of
\cite{SiLeMar} to the singular case, in Corollary
\ref{blockdiagonalREM} and Proposition \ref{blockdiagonalforms}.
 Finally, Section 7 puts
in context our results with related work in the literature. In
particular it is shown how the block-diagonal expression for the
symplectic matrix of Proposition \ref{blockdiagonalforms}
particularizes in the regular case to the normal form obtained in
\cite{SiLeMar}, and a comparison is also made between our results and the
Lagrangian Block-Diagonalization method of Lewis
\cite{Lewisblock}.
\section{Relative Equilibria and simple mechanical
systems}\label{prelim} Let $(\mathcal{P},\w)$ be a smooth finite
dimensional symplectic manifold with symplectic form $\w$ and $G$
a finite dimensional Lie group acting smoothly, properly and in a
Hamiltonian fashion on $(\PP,\w)$ with $\Ad^*$-equivariant
momentum map $\J:\PP\rightarrow\g^*$. Given a $G$-invariant
Hamiltonian function $h\in C^G(\PP)$, a point $z\in\PP$ is called
a \emph{relative equilibrium} (for $h$) if its Hamiltonian
evolution lies inside a group orbit. Equivariance of $\J$ and
Noether's Theorem imply that the Hamiltonian evolution of $z$ is
described as the orbit of $z$ by a one-parameter subgroup of $G$
generated by a Lie algebra element $\xi$ which belongs to
$\g_\mu\subset\g$, where $\mu=\J(z)$ and $\g_\mu$ is the Lie
algebra of $G_\mu$, the stabilizer of $\mu$ for the coadjoint
representation of $G$. The element $\xi$ is called a
\emph{velocity} of the relative equilibrium. Using the usual
notation for the infinitesimal action of $\g$ on $\PP$ the
condition for $z$ to be a relative equilibrium is written as
$X_h(z)=\xi_\PP(z)$, where $X_h$ is the Hamiltonian vector field
associated to the function $h$.

If the stabilizer $G_z$ of $z$ is not discrete then there is a
degeneracy in the choice of a velocity for a given relative
equilibrium, since any representative of the class
$[\xi]\in\g/\g_z$ produces the same orbit of $z$. In any case, by
equivariance of $\J$, the inclusion $G_z\subset G_\mu$ holds.
 The quintuple $(\PP,\w,G,\J,h)$ will be
called in short a symmetric Hamiltonian system.

 The following
definition introduced in \cite{Patrick} is generally adopted as
the correct notion of stability of relative equilibria in
Hamiltonian systems, generalizing in the Hamiltonian context the
concept of Lyapunov stability of fixed equilibria for flows of
vector fields.
\begin{defn}\label{stabdef} A relative
equilibrium $z$ with momentum $\mu=\J(z)$ is said to be
$G_\mu$-stable if for every $G_\mu$-invariant neighbourhood $U$ of
$G_\mu\cdot z$ there exists a neighbourhood $O$ of $z$ such that
the Hamiltonian orbit of $O$ lies in $U$.
\end{defn}
In \cite{LeSi,OrRa99} a method for testing stability of
relative equilibria with singular momentum values in Hamiltonian
systems is developed, generalizing the Energy-Momentum Method of
\cite{Patrick} for relative equilibria with discrete stabilizers.
We quote here the main result, due to its importance in the
subsequent development of the paper. For that, given an element
$\xi\in\g$, define the \emph{augmented Hamiltonian} $h_\xi\in
C^\infty(\PP)$ as \be\label{augmentedhamiltonian}
h_\xi(z)=h(z)-\langle\J(z),\xi\rangle .\ee It well-known that
$z\in\PP$ is a relative equilibrium for the symmetric Hamiltonian
system $(\PP,\w,G,\J,h)$ with velocity $\xi$ if and only if $z$ is
a critical point of $h_\xi$. Also recall that since the
Hamiltonian $G$-action is proper, any stabilizer for this action,
in particular $G_z$, must be compact (see \cite{DuiKol}). Since
$\g_z\subset\g_\mu$ we can choose a $G_z$-equivariant splitting of
$\g_\mu$ as $\g_\mu=\g_z\oplus\g_z^\perp .$
\begin{thm}[\cite{LeSi},\cite{OrRa99}]\label{lermansinger} Let
$(\PP,\w,G,\J,h)$ be a symmetric Hamiltonian system and $z\in\PP$
a relative equilibrium with stabilizer $G_z$ and velocity $\xi$.
Assume that $\J (z)=\mu$ and that $G_\mu$ is compact. Then
$\xi\in\g_\mu$. Let $\xi^\perp$ be the projection of $\xi$ onto
some $G_z$-invariant complement of $\g_z$ in $\g_\mu$ (always
available by compactness of $G_z$). If
${\ed_z^2h_{\xi^\perp}}\rrestr{V_s}$ is definite for some (and
hence any) complement $V_s$ to $\g_\mu\cdot z$ in $\ker T_z\J$,
then $z$ is $G_\mu$-stable.
\end{thm}
In this theorem the ambiguity in the velocity introduced by the
stabilizer of the relative equilibrium appears explicitly. In
typical computations, testing this condition over all possible
$G_z$-invariant complements of $\g_z$ in $\g_\mu$ gives the
sharpest stability results (see \cite{OrRa99} and the example in
Section 5). There is an infinite number of choices for the space
$V_s$ in Theorem \ref{lermansinger}, and any of them is called the
(maximal) \emph{symplectic normal space} at $z$, since it is a
maximal symplectic subspace of the symplectic orthogonal to the
group orbit at $z$.

In this paper we will study a particular case of Hamiltonian
systems of great interest in Classical Mechanics. This is the
class of the so-called symmetric \emph{simple mechanical systems},
which are symmetric Hamiltonian systems where $\PP$ is $T^*Q$, the
cotangent bundle of a smooth, finite-dimensional Riemannian
manifold $(Q,\ll\cdot,\cdot\gg )$ equipped with its canonical
symplectic form $\w$ and $G$ is a finite-dimensional Lie group
acting properly and isometrically on $Q$ and by cotangent lifts on
$T^*Q$.

Following \cite{Smale}, the Hamiltonian function is constructed in
the following way: let $V$ be a smooth $G$-invariant function on
$Q$ and call $\overline{V}=\tau^*V\in C^G(T^*Q)$, where
$\tau:T^*Q\rightarrow Q$ is the cotangent bundle projection. We
will refer to both $\overline{V}$ and $V$ as the \emph{potential
energy}. Let $\FL:TQ\rightarrow T^*Q$ be the \emph{Legendre map}
associated to $\ll\cdot,\cdot\gg$, defined by the formula
\be\label{legendre} \langle\FL (v_x),w_x\rangle=\ll
v_x,w_x\gg,\quad\forall\, v_x,w_x\in T_xQ.\ee The Legendre map is
a $G$-equivariant vector bundle isomorphism covering the identity
on $Q$. With it, we can define the \emph{kinetic energy} $K\in
C^G(\PP)$ as $$ K(p_x)=\frac 12 \ll\FL^ {-1}(p_x),\FL^
{-1}(p_x)\gg,\quad\forall\, p_x\in T_x^*Q .$$ Finally, the
Hamiltonian $h$ is defined by \be\label{smshamiltonian}
h=K+\overline{V}.\ee With respect to the canonical symplectic
structure of $T^*Q$
 the
cotangent-lifted action of $G$ is Hamiltonian, with equivariant
momentum map defined by the expression
\be\label{ctmomentum}\langle\J(p_x),\xi\rangle=\langle
p_x,\xi_Q(x)\rangle,\quad\forall\, \xi\in\g .\ee The symplectic
manifold $T^*Q$ is the phase space of the Hamiltonian system,
while the base $Q$ is usually called \emph{configuration space}.
Accordingly, for any point $p_x$ in $T^*Q$, the projection $x=\tau
(p_x)$ is called the configuration point (or base point) of $p_x$.

A key feature of simple mechanical systems is that both their
geometric and dynamical properties are entirely constructed using
the knowledge of the Riemannian manifold $(Q,\ll\cdot,\cdot\gg)$,
the isometric $G$-action on it and the choice of a potential
energy $V$. Thus, one could reasonably expect that the
implementation on this class of Hamiltonian systems of the
stability test given in Theorem \ref{lermansinger} should simplify
accordingly, and yield easier computations at the level of $Q$ and
$G$ instead of the bigger space $\PP=T^*Q$. The obtention of such
a refinement to simple mechanical systems of this stability test
is the main result of this paper, Theorem \ref{REM}.
\begin{rem}
There are several treatments of this problem in the literature.
Simo \emph{et al.} develop in \cite{SiLeMar} a sophisticated test,
that particularizes Theorem \ref{lermansinger}, called the reduced
Energy-Momentum Method. This test gives sufficient conditions for
the stability of relative equilibria in simple mechanical systems
provided the configuration point of the cotangent re\-la\-ti\-ve
equilibrium has a discrete stabilizer (which is the same as to
require that its momentum value is regular). The main advantage of
this is that it is constructed specifically for this class of
systems, and this fact reflects in less computational difficulties
than the application of the main method, Theorem
\ref{lermansinger}, designed for general Hamiltonian systems. In
\cite{Lewisblock} a Lagrangian analogue of the results of
\cite{SiLeMar} is obtained, being valid also for relative
equilibria of a larger class of mechanical systems. Here, based in
Theorem \ref{lermansinger}, we produce a method for testing
stability of relative equilibria in simple mechanical systems that
could be seen as a generalization of the reduced Energy-Momentum
Method to the singular case, i.e. without requiring discrete
stabilizers of configuration points or regular momentum values.
\end{rem}
\section{A cotangent-bundle adapted splitting of the symplectic
normal space}
 In this section we describe a
realization $V_s\subset T_{p_x}(T^*Q)$ of the symplectic normal
space at a relative equilibrium $p_x$ of a simple mechanical
system, as well as a convenient cotangent-bundle adapted splitting
of $V_s$ which will be extremely useful for the remaining
constructions. Most of the results of this section are merely
expository, and a complete description including proofs and the
obtention of the symplectic normal space at points $p_x$ of
general form, not only relative equilibria, can be found in
\cite{PerRoSD2}.

One of the geometric objects that will be extensively used is the
\emph{locked inertia tensor} $\II$, a family of bilinear positive
semi-definite symmetric forms on $\g$ defined by
\be\label{lockedinertia}
\II(x)(\xi,\eta)=\ll\xi_Q(x),\eta_Q(x)\gg,\quad\forall\,
\xi,\eta\in\g,\,x\in Q.\ee Note that at each point $x$, the kernel
of $\II (x)$ is precisely $\g_x$, the Lie algebra of the
stabilizer of $x$. Therefore $\II(x)$ is a well defined inner
product on $\g$ only at points of $Q$ where the action is locally
free. The locked inertia tensor satisfies the following invariance
and infinitesimal invariance properties (see
\cite{MarLec}):\be\label{propertieslocked1}
 \II(g\cdot
 x)(\Ad_g\xi,\Ad_g\eta)=\II(x)(\xi,\eta)\ee
\be\label{propertieslocked2}
(\D\II\cdot\lambda_Q(x))(\xi,\eta)+\II(x)(\ad_\lambda\xi,\eta)+\II(x)(\xi,\ad_\lambda\eta)=0,\ee
for every $g\in G,\,x\in Q$ and $\xi,\eta,\lambda\in\g$. Note that
$G_{p_x}\subset G_x$ by the equivariance of $\tau$. Let $x\in Q$
be the base point of an element $p_x\in T^*Q$ and denote by
$H=G_x$ its stabilizer. Using the $G$-invariant Riemannian metric
on $Q$ we can  form the splitting \be\label{TQvh} T_xQ= \g\cdot
x\oplus\Sl\ee where $\Sl=(\g\cdot x)^\perp$, and it is usually
called a linear slice (for the $G$-action at $x$). Hence, $\Sl$ is
the space of directions orthogonally complementary to the group
orbit. This is obviously a $H$-invariant splitting for the induced
linear $H$-action on $T_xQ$.

 Next we choose a $G_{p_x}$-invariant splitting of the Lie
algebra $\g$ of the form \be\label{ghr}\g=\h\oplus\lir.\ee This is
always possible by the properness of the $G$-action on $Q$, which
implies that $G_{p_x}$ is compact. A concrete way of choosing
\eqref{ghr} will be introduced in \eqref{gsplitting}. For any
element $\xi\in\g$ we write in a unique way $\xi=\xi^\h+\xi^\lir$,
relative to this splitting. The space $\lir$ collects the elements
of $\g$ that generate nontrivial orbits of $x$.

Noting that $\lir\simeq\g\cdot x$ by the isomorphism
$\xi\mapsto\xi_Q(x)$, we can compose this identification with
\eqref{TQvh} and dualize, to get \be\label{splittangent}
  T_xQ \simeq \lir\oplus \Sl \quad \mathrm{and}\quad
  T_x^*Q  \simeq  \lir^* \oplus \Sl^*.
\ee Associated to the Riemannian structure on $Q$, there is an
Ehresmann connection on $T^*Q$, for which the connection map at
$p_x$, $K:T_{p_x}(T^*Q)\rightarrow T^*_xQ$ is defined as
\be\label{opK}K(X)=\frac{D_c^\nabla}{Dt}\restr{t=0}\hat{c}(t)\quad\forall
X\in T_{p_x}(T^*Q),\ee the covariant differential of $\hat{c}(t)$
along $c(t)=\tau(c(t))$ relative to the Levi-Civita connection
$\nabla$. Here $\hat{c}(t)$ is any local curve
$\hat{c}:(-\epsilon,\epsilon)\rightarrow T^*Q$ projecting to
$c(t)$ and satisfying $i)\,\hat{c}(0)=p_x$, and $ii)\,
\frac{d}{dt}\rrestr{t=0}\hat{c}(t)=X$. This connection map $K$ at
$p_x$ is a $G_{p_x}$-equivariant linear map, which combined with
the differential at $p_x$ of the cotangent bundle projection
$\tau$ yields a $G_{p_x}$-equivariant linear isomorphism
$\Psi:T_{p_x}(T^*Q)\rightarrow T_xQ\oplus T_x^*Q=T^*(T_xQ)$
defined by $$\Psi(X)=(T_{p_x}\tau (X),K(X))\quad\forall\, X\in
T_{p_x}(T^*Q)$$ (see \cite{PerRoSD2}). We call vectors at $p_x$
lying in the kernel of $T_{p_x}\tau$ \emph{vertical}, since they
are tangent to the cotangent fiber through $p_x$. Analogously,
those elements of $T_{p_x}(T^*Q)$ which are in  the kernel of $K$
are called \emph{horizontal}, and are identified through
$T_{p_x}\tau$ with vectors tangent to $Q$ at $x$.  We now compose
the above isomorphism $\Psi$ with the two dual
isomorphisms~\eqref{splittangent} to get a new one
\be\label{sasaki}I:T_{p_x}(T^*Q)\rightarrow (\lir\oplus\Sl)\oplus
(\lir^*\oplus\Sl^*)\ee which can be explicitly expressed as
$I(X)=(\eta,a;\nu,\alpha),$ for the unique $\eta,a,\nu,\alpha$
satisfying
$$\begin{array}{lll}
T_{p_x}\tau (X) & = & a+\eta_Q(x)\\
K(X) & = & \alpha +\FL
\left(\left(\rII^{-1}(\nu)\right)_Q(x)\right).
\end{array}$$
Here $\rII$ denotes the restriction of $\II$ to $\lir$, according
to~\eqref{ghr}. Note that now $\rII$ becomes a well-defined inner
product in $\lir$ and thus also a linear isomorphism
$\rII:\lir\rightarrow\lir^*\simeq \h^\circ$. We can therefore work
in the image of $I$, which we call $I$-representation, instead of
on $T_{p_x}(T^*Q)$, and that is what we will do in the rest of the
paper. Note that in this identification the space of vertical and
horizontal vectors is expressed, respectively, as
$$(0,0;\nu,\alpha)\quad
\forall\,\alpha\in\Sl^*,\,\nu\in\lir^*,\quad\mathrm{and}\quad
 (\eta,a;0,0)\quad\forall\,
a\in\Sl,\,\eta\in\lir .$$ The isomorphism $I$ is
$G_{p_x}$-equivariant with respect to the linear action on the
target space given by \be\label{I-action}g\cdot
(\eta,a;\nu,\alpha)=(\Ad_g\eta,g\cdot
a,\Ad_{g^{-1}}^*\nu,g\cdot\alpha),\ee where $g\cdot a$ and
$g\cdot\alpha$ denote, respectively, the restriction to $G_{p_x}$
of the linear representation of $H$ on $\Sl$ and its
contragredient representation on $\Sl^*$. This action is well
defined since $G_{p_x}\subset H$ and $\Sl$, $\lir$, and their
duals are $H$-invariant by construction.

In order to obtain a convenient characterization of the symplectic
normal space $V_s\subset T_{p_x}(T^*Q)$ at a relative equilibrium
in the $I$-representation, and also for future reference, we quote
some technical results introduced in \cite{PerRoSD2}. There it is
proved that it is possible to extend vectors $v\in\Sl$ to local
vector fields $\overline{v}$ defined in a neighbourhood of $x$, in
a way adapted to the $G$-action,
 and  such  that the family
of vector fields $\lambda_Q,\overline{v}$, for any $\lambda\in\g$
and $v\in\Sl$ spans $T_{x'}Q$ at every $x'$ near $x$.

We will sketch here the obtention of the local field
$\overline{v}$. Recall that by Palais' Tube Theorem \cite{Palais}
we can construct an invariant tubular neighbourhood of an orbit
$G\cdot x\subset Q$ as follows: Let $H=G_x$ act on $G\times \Sl$
as $h\cdot (g,s)=(gh^{-1},h\cdot s)$. Let $G\times_H\Sl$ be the
quotient space for this action. Then there is a $H$-invariant open
ball $U\subset\Sl$ centered at the origin such that the map
$\sigma:G\times_HU\rightarrow Q$ defined by \be\label{tube}\sigma
([g,s])=g\cdot\exp_x s\ee is a diffeomorphism onto a $G$-invariant
neighbourhood $O$ of $G\cdot x$. Here $\exp_x:T_xQ\rightarrow Q$
denotes the exponential map associated to $\ll\cdot,\cdot\gg$.
This diffeomorphism is $G$-equivariant with respect to the given
action on $Q$ and the $G$-action on $G\times_H \Sl$ defined by
$g'\cdot [g,s]=[g'g,s]$.

 Now choose any inner product on $\g$ such that the splitting
 \eqref{ghr} is an orthogonal direct sum. Extend this inner product
 by right translations to a $H$-invariant Riemannian metric on $G$. Then
we can interpret $\lir$ as a linear slice at the identity for the
free $H$-action on $G$ given by $(h,g)\mapsto gh^{-1}$. It follows
that if we call $\exp_e$ the exponential map for this metric on
$G$, there is a small $H$-invariant neighbourhood of $e$ in $G$
such that every $g$ belonging to it can be written as
$g=\exp_e\xi^\lir h^{-1}$ for unique elements $\xi^\lir\in\lir$
and $h\in H$. Using the Tube Theorem every element in a small
neighbourhood of $x$ in $Q$ (not $G$-invariant in general) can be
expressed as $x'=\sigma([\exp_e\xi^\lir h^{-1},s])$ for unique
elements $\xi^\lir\in\lir,\, h\in H$ and $s\in U$. It is easy to
prove that for any $x'=\sigma ([\exp_e\xi^\lir h^{-1},s])$ near
$x$, the formula $$ F_v^t(x')=\sigma ([\exp_e\xi^\lir
h^{-1},s+th\cdot v])$$ defines a flow $F^t_v$ for any $v\in
U\subset\Sl$. The associated local field $\overline{v}$ is then defined
as
\be\label{localv}\overline{v}(x')=\frac{d}{dt}\restr{t=0}F_v^t(x')\ee
for every $x'$ near $x$. In the following theorem we collect the
most important properties of this family of local vector fields.
For the proof, see \cite{PerRoSD2}.
\begin{thm}\label{covariant} Let $v,w\in\Sl$ and
$\eta,\lambda,\xi,\xi_i,\xi_j\in\g$. Then
\begin{itemize}
\item[(i)] $ \ll \nabla_{{\xi_i}_Q}{\xi_j}_Q(x),\lambda_Q(x)\gg =
\frac 12\left\{ \left(\mathbf{D}\II \cdot
{\xi^\lir_i}_Q(x)\right)(\xi_j,\lambda)-\II
(x)(\xi^\lir_i,\left[\xi_j,\lambda\right])\right\}$

\item[(ii)] $\ll \nabla_{{\xi_i}_Q}{\xi_j}_Q(x),w\gg=-\frac
12\left(\mathbf{D}\II\cdot w\right)(\xi_i^\lir,\,\xi_j)$

\item[(iii)] $\ll
\nabla_{{\xi}_Q}\,\overline{v}(x),\lambda_Q(x)\gg=\frac
12\left(\mathbf{D}\II\cdot v\right)(\xi^\lir ,\,\lambda)$

\item[(iv)] $\ll
\nabla_{\overline{v}}\,{\xi}_Q(x),\lambda_Q(x)\gg=\frac
12\left(\mathbf{D}\II \cdot v\right)(\xi,\lambda)$

\item[(v)] $\ll
\nabla_{\overline{v}}\,{\xi}_Q(x),w\gg=\ll\nabla_{{\xi}_Q}\,\overline{v}(x),w\gg+\ll\xi^\h\cdot
v,w\gg $.
\end{itemize}
Here, $\xi^\lir,\xi_i^\lir$ and $\xi^\h$ denote the projections of
elements of $\g$ onto $\lir$ and $\h$ according to \eqref{ghr}.
\end{thm}
\paragraph{{\bf Notation:}} We will introduce for any $v\in\Sl$ a linear map
$C(v):\lir\rightarrow\Sl^*$ defined as \be\label{defC}\ll
C(v)(\xi^\lir ),w\gg_\Sl=\ll\nabla_{\xi_Q}\,\overline{v}(x),w\gg
,\ee where $\ll\cdot,\cdot\gg_\Sl$ is the restriction of the
metric on $Q$ to $\Sl\subset T_xQ$. Note that $C$ is not linear in
$v$ since it depends on the concrete extension $\overline{v}$. We
will also employ the following notation: if $W$ is a linear
subspace of the linear space $V$ and $\iota:W\hookrightarrow V$
its inclusion, we will write $\Proj_W:V^*\rightarrow W^*$ for its
dual projection.

We fix from now on a point of the form $p_x=\FL(\xi_Q(x))$ with
$G_x=H$. It follows from \eqref{ctmomentum} and the definition of
the locked inertia tensor that $p_x$ has momentum $\mu=\J
(p_x)=\II(x)(\xi)=\rII(\xi^\lir)$. The reason for this choice of
$p_x$ will be clear in the following section, where it is
explained why every relative equilibrium of a simple mechanical
system must be precisely of this form.
\paragraph*{\bf{Remark.}}
It is a well-known property of points of the form
$p_x=\FL(\xi_Q(x))$ that $G_{p_x}=H\cap G_\mu$. This follows
immediately from the relation $G_{p_x}=H_{p_x}$ and identifying
$p_x=(0,\mu)\in \Sl^*\oplus\lir^*$ using the $H$-isomorphism
\eqref{splittangent}. For general points of $T^*Q$ one has only an
inclusion $G_{p_x}\subset H\cap G_\mu$.

We now make a concrete choice for the complement $\lir$ in
\eqref{ghr}, as well as for other relevant subspaces of $\g$.
Start by choosing a $G_{p_x}$-invariant complement $\mathfrak{p}$
to $\g_{p_x}$ in $\g_\mu$, i.e.
\begin{equation}\label{optimal}\g_\mu=\g_{p_x}\oplus\mathfrak{p}.\end{equation}
Next, define a $G_{p_x}$-invariant complement $\mathfrak{t}$ to
$\h\oplus\mathfrak{p}$ in $\g$ in such a way that defining
$\lir=\mathfrak{p}\oplus\mathfrak{t}$ we have that $\mathfrak{p}$
and $\mathfrak{t}$ are orthogonal with respect to the restricted
locked inertia tensor $\rII$. We can then write
\begin{equation}\label{gsplitting}\g=\h\oplus\lir=\h\oplus\mathfrak{p}\oplus\mathfrak{t},\end{equation} and
hence we have constructed the splitting $\eqref{ghr}$.

Let us define the following subspace of $\g$ \be\label{qmu}
\q^\mu=\{\lambda\in\mathfrak{t}\, : \,\Proj_\h\, [\ad^*_\lambda\mu
]=0 \} .\ee This space will play an important role in our
characterization of $V_s$, and it can be proved (see
\cite{PerRoSD2}) that it is isomorphic to the symplectic normal space at
$\mu$ for the restriction to $H$ of the coadjoint action of $G$ on
$\OO_\mu$, the coadjoint orbit containing $\mu$. Note also that by
using \eqref{propertieslocked2} we can write
$\q^\mu=\{\lambda\in\mathfrak{t}\, : \,\Proj_\h\, [(\D\II\cdot
\lambda_Q(x))(\xi^\lir)=0]\}$.

As a particular case of Theorem 6.1 in \cite{PerRoSD2}  the space
$\ker T_{p_x}\J$ consists in the $I$-representation in the
elements $(\eta,a;\nu,\alpha)\in (\lir\oplus\Sl)\oplus
(\lir^*\oplus\Sl^*)$ satisfying
$$\langle\nu,\lambda^\lir\rangle-\frac
12\left\{\left(\D\II\cdot\xi_Q(x)\right)(\lambda^\lir,\eta)-\left(\D\II\cdot
a\right)(\lambda^\lir,\xi)-\langle\ad^*_{\lambda^\lir}\mu,\eta\rangle\right\}
+\langle\ad^*_{\lambda^\h}\mu,\eta\rangle=0$$ for every
$\lambda^\lir\in\lir$ and $\lambda^\h\in\h$. Also, in the
$I$-representation
$$\g_\mu\cdot p_x=\left\{\left(\lambda,0;\frac
12\Proj_\lir\left[(\D\II\cdot\xi_Q(x))(\lambda)\right],-\frac
12\Proj_\Sl\left[(\D\II\cdot
(\cdot))(\lambda,\xi)\right]\right)\,:\,\forall\,\lambda\in\mathfrak{p}
\right\}.$$ From the above two expressions it is easy to obtain
that the symplectic normal space $V_s$, a complement to
$\g_\mu\cdot p_x$ in $\ker T_{p_x}\J$, can be chosen to be
\be\label{Vre}\begin{array}{lll}V_s  & = & \left\{ \left(\lambda,
a ; \frac 12\Proj_\lir\, \left[(\D\II
\cdot\xi^\lir_Q(x))(\lambda)+\ad^*_\lambda\mu-(\D\II\cdot
a)(\xi^\lir)\right],\gamma\right.\right.\vspace{5pt}\\
 & - & \left.\left.\frac 12 \Proj_\Sl\,\left[(\D\II\cdot (\cdot
))(\xi^\lir,\lambda)\right]+\ll C(a)(\xi^\lir),\cdot\gg_\Sl
\right)\,:\,\forall\,\lambda\in\q^\mu,\,a\in\Sl,\,\gamma\in\Sl^*\right\},
\end{array}\ee
with $\q^\mu$ defined in \eqref{qmu} and $\Sl$ in \eqref{TQvh}.
The symplectic normal space $V_s$ is $G_{p_x}$-invariant by
construction with respect to the action \eqref{I-action} in its
ambient space (see \cite{PerRoSD2}).
\section{Stability of singular relative equilibria in simple
mechanical systems} In the following we will be in the setup of
Section \ref{prelim} and fix a simple mechanical system
$h=K+\overline{V}$ as in \eqref{smshamiltonian}. In this
framework, once an element $\xi\in\g$ is chosen, we can separate
the augmented Hamiltonian \eqref{augmentedhamiltonian} into a
kinetic and a potential part as
$$h_\xi=K_\xi+\overline{V_\xi},\quad \text{where}$$
$$ K_\xi
(p_x) =  \frac 12 \| p_x-\chi^\xi(x) \|^2\quad\text{and}\quad
V_\xi(x)  =  V(x)-\frac 12 \II(x)(\xi,\xi).$$ As for
the potential energy we have used the notation
$\overline{V_\xi}=\tau^*V_\xi$. The one-form $\chi^\xi$ is defined
by \be \label{chi}\chi^\xi(x)=\FL(\xi_Q(x)).\ee The functions
$K_\xi$ and $V_\xi$ are called the \emph{augmented kinetic energy}
and \emph{augmented potential energy} respectively. Recall now
(see for instance Theorem 4.1.2 and Proposition 4.2.1 in
\cite{MarLec}) that with the introduction of these two auxiliary
functions we have the following characterization of relative
equilibria:

\begin{prop}\label{charactre} Let $p_x\in
T^*Q$. The following are equivalent:\begin{itemize} \item[(i)]
$p_x$ is a relative equilibrium of \eqref{smshamiltonian} with
momentum $\mu$ and velocity $\xi\in \g_\mu$ \item[(ii)] $p_x$ is a
critical point of the augmented Hamiltonian $h_\xi$ \item[(iii)]
$p_x$ is simultaneously a critical point of $K_\xi$ and
$\overline{V_\xi}$ \item[(iv)] $p_x=\FL(\xi_Q(x))$ and $x$ is a
critical point of $V_\xi$.
\end{itemize}
\end{prop}
Note that (iv) restricts the form of phase space points candidates
to be relative equilibria of \eqref{smshamiltonian} to be of the
form $p_x=\FL(\xi_Q(x))$. Thus is the reason for studying in detail in the previous section the symplectic normal space only at this class
of points.

Let $p_x=\FL(\xi_Q(x))\in T^*_xQ$ be a relative equilibrium for
the simple mechanical system \eqref{smshamiltonian} with momentum
$\J(p_x)=\mu$. We call $H=G_x$ and we choose a $(G_{p_x}=H\cap
G_\mu)$-invariant inner product on $\g$ relative to which we
construct the splittings \eqref{optimal} and \eqref{gsplitting}.
 Note that by hypothesis $x$ is a critical point
of $V_{\xi'}$ for any $\xi'\in\g_\mu$ such that
$[\xi'-\xi]=0\in\g_\mu/\g_{p_x}$. In other words, any such $\xi'$
is a velocity for the relative equilibrium $p_x$. In particular,
this happens for $\xi^\perp\in\mathfrak{p}$, the projection of
$\xi$ onto $\mathfrak{p}$ according to \eqref{optimal}.

Let $\delta p\in T_{p_x}(T^*Q)$ be a tangent vector at $p_x$. We
will write its horizontal and vertical components as $\delta
p^H=T_{p_x}\tau (\delta p)\in T_xQ$ and $\delta p^V=K(\delta p)\in
T_x^*Q$ respectively. It is clear that if $I(\delta
p)=(\lambda,b;\nu,\gamma)$ is the $I$-representation of $\delta p$
one has that $$\begin{array}{lll}
\delta p^H & = & \lambda_Q(x)+b\\
\delta p^V & = & \FL\left((\rII^{-1}(\nu))_Q(x)\right)+\gamma.
\end{array}$$
Also, using \eqref{splittangent} we can express the horizontal and
vertical variations $\delta p^H$ and $\delta p^V$ as elements of
$\lir\oplus\Sl$ and $\lir^*\oplus\Sl^*$ respectively like
$$\begin{array}{lll}
\delta p^H & = & (\lambda,b)\\
\delta p^V & = & (\nu,\gamma).
\end{array}$$
We will use both notations indistinctly.

Finally, for a curve $c(t)\in Q$ with $c(0)=x$ we will write
$\Hor_{p_x} (c(t))$ for its horizontal lift to $T^*Q$ at the point
$p_x$ with respect to the Levi-Civita connection. Equivalently,
$\Hor_{p_x} (c(t))$ is the parallel translation of $p_x$ along the
curve $c(t)$.
\begin{lemma}\label{hessians}Let $p_x=\FL(\xi_Q(x))\in T^*_xQ$ be a relative equilibrium
for the simple mechanical system \eqref{smshamiltonian} with
momentum $\mu$ and velocity $\xi\in\g_\mu$. Let $\xi^\perp$ the
orthogonal projection of $\xi$ onto $\mathfrak{p}$ according to
\eqref{optimal}. Then, for any $\delta p_1,\delta p_2\in
T_{p_x}(T^*Q)$
\begin{itemize}
\item[(i)] $\ed^2_{p_x} \overline{V_{\xi^\perp}}(\delta p_1,\delta
p_2) = \ed^2_x V_{\xi^\perp}(\delta p_1^H,\delta p_2^H)$
\item[(ii)] $\ed^2_{p_x} K_{\xi^\perp} (\delta p_1,\delta p_2) =
\ll \delta p_1^V-(T_x\chi^{\xi^\perp}\cdot\delta p_1^H)^V,\delta
p_2^V-(T_x\chi^{\xi^\perp}\cdot\delta p_2^H)^V\gg $.
\end{itemize}
\end{lemma}
\begin{proof}
(i) follows immediately by the definition of $\delta p^H$ and
noting that $\overline{V_{\xi^\perp}}=\tau^* V_{\xi^\perp}$. To
prove (ii) we will consider horizontal and vertical vectors
separately. Let $I(\delta p_i)=(0;\delta p_i^V)$ for $i=1,2$. Then
$$\begin{array}{lll} \ed^2_{p_x} K_{\xi^\perp}
(\delta p_1,\delta p_2) & = & \frac 12\frac
{d}{ds}\rrestr{s=0}\frac {d}{dt}\rrestr{t=0} \parallel p_x+t\delta
p_1^V+s\delta p_2^V-\chi^{\xi^\perp} (x)\parallel^2\vspace{1mm}\\
 & = & \frac
12\frac {d}{ds}\rrestr{s=0}\frac {d}{dt}\rrestr{t=0}
\parallel t\delta p_1^V+s\delta p_2^V
(x)\parallel^2\vspace{1mm}\\
& = & \ll\delta p_1^V,\delta p_2^V\gg .
\end{array}$$
If $I(\delta p_1)=(0;\delta p_1^V)$ and $I(\delta p_2)=(\delta
p_2^H;0)$, let $c^t_{\delta p_2^H}(x)$ be any smooth curve
satisfying $c^0_{\delta p_2^H}(x)=x$ and
$\frac{d}{dt}\rrestr{t=0}c^t_{\delta p_2^H}(x)=\delta p_2^H$. Then
$$\begin{array}{lll} \ed^2_{p_x} K_{\xi^\perp}
(\delta p_1,\delta p_2) & = & \frac 12\frac
{d}{ds}\rrestr{s=0}\frac {d}{dt}\rrestr{t=0}
\parallel \Hor_{p_x+t\delta p_1^V}(c^s_{\delta p_2^H}(x))-\chi^{\xi^\perp}(c^s_{\delta
p_2^H}(x))\parallel^2\\
 & = & -\frac {d}{dt}\rrestr{t=0} \ll \nabla_{\delta p_2^H}\chi^{\xi^\perp}
  (x), p_x+t\delta p_1^V-\chi^{\xi^\perp}(x) \gg .\end{array}$$
But since $p_x$ is a relative equilibrium with velocity $\xi$,
then $p_x= \chi^{\xi^\perp} (x)$ and then
$$\ed^2_{p_x} K_{\xi^\perp}
(\delta p_1,\delta p_2)= -\ll \nabla_{\delta
p_2^H}\chi^{\xi^\perp} (x),\delta p_1^V\gg .
$$
Finally, consider two variations of the form $I(\delta
p_1)=(\delta p_1^H;0)$ and $I(\delta p_2)=(\delta p_2^H;0)$. Then
$$\begin{array}{l}
\ed^2_{p_x} K_{\xi^\perp} (\delta p_1,\delta p_2)=\vspace{1mm}\\
  = \frac
12\frac {d}{ds}\rrestr{s=0}\frac{d}{dt}\rrestr{t=0}
\parallel \Hor_{\Hor_{p_x}(c^t_{\delta p_1^H}(x))}(F^s_{\widetilde{\delta
p_2^H}}(c^t_{\delta p_1^H}(x)))-\chi^{\xi^\perp}
(F^s_{\widetilde{\delta
p_2^H}}(c^t_{\delta p_1^H}(x)))\parallel^2\\
 =  -\frac{d}{dt}\rrestr{t=0}\ll \nabla_{\widetilde{\delta
p_2^H}}\chi^{\xi^\perp} (c^t_{\delta
p_1^H}(x)),\Hor_{p_x}(c^t_{\delta
p_1^H}(x))-\chi^{\xi^\perp}(c^t_{\delta p_1^H}(x))\gg\\
 = \ll
\nabla_{\delta p_2^H}\chi^{\xi^\perp} (x),\nabla_{\delta
p_1^H}\chi^{\xi^\perp} (x)\gg .
\end{array}$$
Where if $\delta p_2=(\lambda,b)$, then $\widetilde{\delta p_2^H}$
is the local vector field $\overline{b}+\lambda_Q$, and
$F^s_{\widetilde{\delta p_2^H}}$ denotes its flow.
 Recalling now the definition \eqref{opK} of
the operator $K$, for any $v\in T_xQ$ and local extension
$\widetilde{v}$ one has
$$\nabla_{v}\,\chi^{\xi^\perp} (x)=\nabla_{\widetilde{v}}\,\chi^{\xi^\perp} (x)=K(T_x\chi^{\xi^\perp}\cdot v)=(T_x\chi^{\xi^\perp}\cdot v)^V ,$$
and the result is proved.
\end{proof}
Since every horizontal variation $\delta p^H$ can be written as
the sum of two contributions one coming from $\lir$ and the other
from $\Sl$, that is $\delta p^H=\zeta_Q(x) + b$, for
$\zeta\in\lir,\, b\in\Sl$, we can consider these two contributions
separately and thus obtain concrete expressions for
$(T_x\chi^{\xi^\perp}\cdot\delta p^H)^V$.
\begin{lemma}\label{hessians2}
Let $\zeta\in\lir$ and $b\in\Sl$, and identify $T_x^*Q$ with
$\lir^*\oplus\Sl^*$ by the isomorphism~\eqref{splittangent}. Then,
\begin{itemize}
\item[(i)] $\begin{array}{l}(T_x\chi^{\xi^\perp}\cdot\zeta_Q(x))^V
=\vspace{1mm}
\\ \left(\frac 12 \Proj_\lir\,\left[\left(\D\II\cdot
\xi^\perp_Q(x)\right)(\zeta)+\ad^*_\zeta\mu+2\left(\D\II\cdot\zeta_Q(x)\right)(\xi^\perp)\right],
-\frac 12\Proj_\Sl\,\left[\left(\D\II\cdot
(\cdot)\right)(\xi^\perp,\zeta)\right]\right).\end{array}$\vspace{1mm}
\item[(ii)] $(T_x\chi^{\xi^\perp}\cdot b)^V  =\left(\frac 12
\Proj_\lir\,\left[\left(\D\II\cdot b\right)(\xi^\perp)\right],\ll
C(b)(\xi^\perp),\cdot\gg_\Sl \right)$.
\end{itemize}
\end{lemma}
\begin{proof}The proof is just an immediate consequence of
Theorem \ref{covariant}, in particular of items (i), (ii), (iv),
and (v). We will just prove (ii). Recalling that $\nabla$ is a
metric connection, then
$$(T_x\chi^{\xi^\perp}\cdot
b)^V=\nabla_{\overline{b}}\left(\FL(\xi^\perp_Q)\right)(x)=\FL(\nabla_{\overline{b}}\,\xi^\perp_Q)(x).$$
For $\lambda\in\lir$ we have, from item (iv) of Theorem
\ref{covariant}
$$\langle (T_x\chi^{\xi^\perp}\cdot
b)^V,\lambda_Q(x)\rangle
=\ll\nabla_{\overline{b}}\,\xi^\perp_Q(x),\lambda_Q(x)\gg =\frac
12(\D\II\cdot b)(\xi^\perp,\lambda).$$  Similarly, if $w\in\Sl$,
then by item (v) of the above referred theorem we have
$$\langle (T_x\chi^{\xi^\perp}\cdot
b)^V,w\rangle =\ll\nabla_{\overline{b}}\,\xi^\perp_Q(x),w\gg =\ll
C(b)(\xi^\perp),w\gg_\Sl .$$ This yields result (ii) The proof of
(i) is identical, including some manipulations using the
infinitesimal equivariance property of the locked inertia tensor
given in \eqref{propertieslocked2}.
\end{proof}
 We apply now the
results obtained so far in order to produce a singular version of
the reduced Energy-Momentum Method of \cite{SiLeMar}. Consider $\PP=T^*Q$ with its
canonical symplectic form in
the statement of Theorem \ref{lermansinger}, where $(Q,\ll\cdot,\cdot\gg )$ is a
Riemannian manifold on which the Lie group $G$ acts by isometries.
The Hamiltonian action on $T^*Q$ is the cotangent lift of the
action on $Q$, and the Hamiltonian is given by
\eqref{smshamiltonian} defining a simple mechanical system. Let
$\xi$ be an element of the Lie algebra $\g$ of $G$. Fix $p_x=\FL
(\xi_Q(x))$ with momentum $\mu$. Assume that
\begin{itemize} \item[(i)] $x$ is a critical point of $V_\xi$, and
\item[(ii)] $G_\mu$ is compact. \end{itemize} Then $p_x$ is a
relative equilibrium for our simple mechanical system with
velocity $\xi$ satisfying $\xi\in\g_\mu$ and Theorem
\ref{lermansinger} can be applied to study its $G_\mu$-stability.
That is, we need to determine when
$\ed^2_{p_x}h_{\xi^\perp}\rrestr{V_s}$ is definite. For that, we
will use the characterization of the symplectic normal space
$V_s$, given in \eqref{Vre} which establishes a linear isomorphism
$\kappa:\q^\mu\oplus \Sl\oplus \Sl^*\rightarrow V_s\subset
T_{p_x}T^*Q\simeq (\lir\oplus\Sl )\oplus (\lir^*\oplus\Sl^*)$
explicitly expressed as
\be\label{isore2}\begin{array}{lll}\kappa\,(\lambda,a,\gamma)=
\left(\lambda, a \right.& ; &\left. \frac 12\Proj_\lir\,
\left[(\D\II
\cdot\xi^\perp_Q(x))(\lambda)+\ad^*_\lambda\mu-(\D\II\cdot
a)(\xi^\perp)\right],\right.\vspace{2mm}\\ &  & \left.\gamma
-\frac 12 \Proj_\Sl\,\left[(\D\II\cdot (\cdot
))(\xi^\perp,\lambda)\right]+\ll C(a)(\xi^\perp),\cdot\gg_\Sl
\right).
\end{array}\ee
This map is $G_{p_x}$-equivariant with respect to the action on
$\q^\mu\oplus\Sl\oplus\Sl^*$ given by $$g\cdot
(\lambda,a,\gamma)=(\Ad_g\lambda,g\cdot a,g\cdot \gamma),$$ and
the action \eqref{I-action} on $(\lir\oplus\Sl)\oplus
(\lir^*\oplus\Sl^*)$.
\begin{defn}\label{correction}
At a relative equilibrium $p_x=\FL(\xi_Q(x))$, with momentum $\mu$
velocity $\xi\in\g_\mu$, define the linear subspace $\Sigma$ of
$T_xQ$ isomorphic to $\q^\mu\oplus \Sl$, as \be\label{sigmadef}
\Sigma=\left\{\lambda_Q(x)+a\in T_xQ\, :\, \lambda\in\q^\mu ,\,
a\in\Sl\right\} ,\ee where $\q^\mu$ is defined in \eqref{qmu}.
\\
For any two vectors $v_1,v_2\in T_xQ$, define the correction term
as the symmetric bilinear form on $T_xQ$ defined by
\be\label{correctionterm} \corr_\xi (x)
(v_1,v_2)=\langle\Proj_\lir\left[(\D\II\cdot v_1)(\xi
)\right],\rII^{-1}\left(\Proj_\lir\left[(\D\II\cdot v_2)(\xi
)\right]\right) \rangle .\ee
\end{defn}
We can now state the main result of this section.
\begin{thm}[reduced Energy-Momentum Method]\label{REM}Suppose that $p_x=\FL(\xi_Q(x))$ is a
relative equilibrium of the simple mechanical system
\eqref{smshamiltonian} with momentum $\mu$ and velocity
$\xi\in\g_\mu$. Assume that $G_\mu$ is compact, $G_x=H$ and that a
$(H\cap G_\mu)$-invariant splitting
$\g_\mu=\g_{p_x}\oplus\mathfrak{p}$. Let $\xi^\perp$ be the
orthogonal projection of $\xi$ onto $\mathfrak{p}$. Then if
\begin{itemize} \item[(i)] $\dim Q-\dim G+\dim G_x >0$, and
\item[(ii)] $\left(\ed^2_{x}V_{\xi^\perp}+\corr_{\xi^\perp}
(x)\right)\rrestr{\Sigma}$ is positive definite
\end{itemize}
or
\begin{itemize} \item[(i)] $\dim Q-\dim
G+\dim G_x =0$, and \item[(ii)]
$\left(\ed^2_{x}V_{\xi^\perp}+\corr_{\xi^\perp}
(x)\right)\rrestr{\Sigma}$ is definite (positive or negative)
\end{itemize}
then the relative equilibrium is $G_\mu$-stable.
\end{thm}
\begin{proof} According to Theorem \ref{lermansinger} the relative equilibrium is stable provided
${\ed_z^2h_{\xi^\perp}}\rrestr{V_s}$ is definite. As both
$K_{\xi^\perp}$ and $\overline{V_{\xi^\perp}}$ have a critical
point at $p_x$ ($(iii)$ in Proposition \ref{charactre}), and since
$h_{\xi^\perp}=K_{\xi^\perp}+\overline{V_{\xi^\perp}}$, then
$${\ed_{p_x}^2h_{\xi^\perp}}\rrestr{V_s}=\left(\ed_{p_x}^2K_{\xi^\perp}+\ed_{p_x}^2\overline{V_{\xi^\perp}}\right)\rrestr{V_s}.$$
 With the
isomorphism $\kappa$ in \eqref{isore2} we can compute each of
these Hessians in $V_s$ parameterized by elements in
$\q^\mu\oplus\Sl\oplus\Sl^*$.

Let us compute first $\ed_{p_x} K_{\xi^\perp}\rrestr{V_s}$. Let
$\lambda,\lambda_1,\lambda_2\in  \q^\mu,\,a,a_1,a_2\in\Sl$ and
$\beta,\beta_1,\beta_2\in\Sl^*$. Note that using \eqref{isore2} we
have that the vertical variation $\delta p^V$ corresponding to an
element $\lambda\in\q^\mu$, i.e. $I(\delta
p)=\kappa\,(\lambda,0,0)$, is
$$\delta p^V=\left(\frac
12\Proj_\lir\, \left[(\D\II\cdot
(\xi^\perp)_Q(x))(\lambda)+\ad^*_\lambda\mu\right],-\frac
12\Proj_\Sl [(\D\II\cdot (\cdot ))(\xi^\perp,\lambda)]\right).$$
Using Lemma \ref{hessians2} we have
$$\delta p^V-(T_x\chi^{\xi^\perp}\cdot\delta p^H)^V=\left(-\Proj_\lir\, \left[(\D\II\cdot \lambda_Q(x))(\xi^\perp)\right],0\right).$$
Similarly,  if $a\in\Sl$ we have for $I(\delta p)=\kappa\,(0,a,0)$
$$\delta p^V=\left(-\frac 12\Proj_\lir\, \left[(\D\II\cdot a)(\xi^\perp)\right],\ll C(a)(\xi^\perp),\cdot\gg_\Sl\right),$$
and then
$$\delta p^V-(T_x\chi^{\xi^\perp}\cdot\delta p^H)^V=\left(-
\Proj_\lir\, \left[(\D\II\cdot a)(\xi^\perp)\right],0\right).$$
Finally, if $\beta\in\Sl^*$ for $I(\delta p)=\kappa\,(0,0,\beta)$
we have $\delta p^V=(0,\beta)$ and $\delta p^H=(0,0)$ so
$$\delta p^V-(T_x\chi^{\xi^\perp}\cdot\delta p^H)^V=(0,\beta).$$

 According to (ii) in Lemma
\ref{hessians} we can now write, for
$\lambda,\lambda_1,\lambda_2\in\q^\mu,\, a,a_1,a_2\in\Sl$ and
$\beta,\beta_1,\beta_2\in\Sl^*$
$$\begin{array}{lll}
\ed_{p_x}^2K_{\xi^\perp}
\left(\kappa(\lambda_1,0,0),\kappa(\lambda_2,0,0)\right)& = &
\langle \Proj_\lir [(\D\II\cdot
{\lambda_1}_Q(x))(\xi^\perp)],\rII^{-1} \Proj_\lir
[(\D\II\cdot {\lambda_2}_Q(x))(\xi^\perp)] \rangle\vspace{1mm}\\
\ed_{p_x}^2K_{\xi^\perp} (\kappa(\lambda,0,0),\kappa(0,a,0))& = &
\langle \Proj_\lir [(\D\II\cdot
\lambda_Q(x))(\xi^\perp)],\rII^{-1} \Proj_\lir
[(\D\II\cdot a)(\xi^\perp)] \rangle\vspace{1mm}\\
\ed_{p_x}^2K_{\xi^\perp} (\kappa(0,a_1,0),\kappa(0,a_2,0))& = &
\langle \Proj_\lir [(\D\II\cdot a_1)(\xi^\perp)],\rII^{-1}
\Proj_\lir
[(\D\II\cdot a_2)(\xi^\perp)] \rangle\vspace{1mm}\\
\ed_{p_x}^2K_{\xi^\perp}
(\kappa(0,0,\beta_1),\kappa(0,0,\beta_2))& = & \ll
\beta_1,\beta_1\gg_{\Sl^*}\vspace{1mm}\\
\ed_{p_x}^2K_{\xi^\perp} (\kappa(0,0,\beta),\kappa(\lambda,0,0))&
= & 0\vspace{1mm}\\
\ed_{p_x}^2K_{\xi^\perp} (\kappa(0,0,\beta),\kappa(0,a,0))& = &
0.
\end{array}$$
Where $\ll\cdot,\cdot\gg_{\Sl^*}$ is the inner product in $\Sl^*$
induced from $\ll\cdot,\cdot\gg_\Sl$ via the Riemannian Legendre
map \eqref{legendre}. We compute now the remaining contribution,
the Hessian of the augmented potential energy. Using Lemma
\ref{hessians} it is immediate to obtain
$$\begin{array}{lll}
\ed^2_{p_x}\overline{V_{\xi^\perp}}(\kappa
(\lambda_1,a_1,0),\kappa (\lambda_2,a_2,0)) & = &
\ed^2_xV_{\xi^\perp} ({\lambda_1}_Q(x)+a_1,{\lambda_2}_Q(x)+a_2)\\
\ed^2_{p_x}\overline{V_{\xi^\perp}}(\kappa (\lambda,a,0),\kappa (0,0,\beta)) & = & 0\\
\ed^2_{p_x}\overline{V_{\xi^\perp}}(\kappa(0,0,\beta_1),\kappa(0,0,\beta_2))
& = & 0
\end{array}$$
for every $\lambda,\lambda_1\lambda_2\in\q^\mu,\,
a,a_1,a_2\in\Sl,\, \beta,\beta_1,\beta_2\in\Sl^*$.

Therefore, ${\ed_{p_x}^2h_{\xi^\perp}}\rrestr{V_s}$
block-diagonalizes in the two blocks
$\kappa(\q^\mu\oplus\Sl\oplus\{0\})$ and $\kappa
(\{0\}\oplus\{0\}\oplus\Sl^*)$ as
$$\begin{array}{lll}
{\ed_{p_x}^2h_{\xi^\perp}}\rrestr{V_s}(\kappa
(\lambda_1,a_1,0),\kappa (\lambda_2,a_2,0)) & =
& \left(\ed^2_xV_{\xi^\perp}+\corr_{\xi^\perp} (x)\right) ({\lambda_1}_Q(x)+a_1,{\lambda_2}_Q(x)+a_2)\\
{\ed_{p_x}^2h_{\xi^\perp}}\rrestr{V_s}(\kappa (\lambda,a,0),\kappa (0,0,\beta)) & = & 0\\
{\ed_{p_x}^2h_{\xi^\perp}}\rrestr{V_s}(\kappa (0,0,\beta_1),\kappa
(0,0,\beta_2)) & = & \ll \beta_1,\beta_1\gg_{\Sl^*}
\end{array}$$
for every $\lambda,\lambda_1\lambda_2\in\q^\mu,\,
a,a_1,a_2\in\Sl,\, \beta,\beta_1,\beta_2\in\Sl^*$. From the above
expression and Definition \ref{correction} it follows that the
bilinear form ${\ed_{p_x}^2h_{\xi^\perp}}\rrestr{V_s}$ is
equivalent to the pair of bilinear forms
$\left(\ed^2_xV_{\xi^\perp}+\corr_{\xi^\perp}
(x)\right)\rrestr{\Sigma}$ and $\ll\cdot,\cdot\gg_{\Sl^*}$.

 Let us examine now when ${\ed_{p_x}^2h_{\xi^\perp}}\rrestr{V_s}$ is definite. Since the block
$\ll\cdot,\cdot\gg_{\Sl^*}$ is positive-definite or trivial, there
are two possible scenarios. (i) $\dim \Sl= 0$ and (ii) $\dim \Sl>
0$. Note that $\dim \Sl=0$ if and only if the dimension of the
orbit $G\cdot x$ equals $\dim Q$. Recall also that $\dim G\cdot x
=\dim G-\dim G_x$. In this case
$\ed^2_{p_x}h_{\xi^\perp}\rrestr{V_s}$ consists only in the block
$\left(\ed^2_xV_{\xi^\perp}+\corr_{\xi^\perp}(x)\right)\rrestr{\Sigma}$
and the relative equilibrium is $G_\mu$-stable
 provided this block is definite (positive or negative).

For the other possibility, if $\dim \Sl\neq 0$ then $\dim Q
>\dim G\cdot x$ and the block $\ll\cdot,\cdot\gg_{\Sl^*}$ is positive definite since the metric on $Q$ is Riemannian, so
${\ed_{p_x}^2h_{\xi^\perp}}\rrestr{V_s}$ is definite if and only
if the block $\left(\ed^2_xV_{\xi^\perp}+\corr_{\xi^\perp}
(x)\right)\rrestr{\Sigma}$ is positive definite. This completes
the proof of the theorem.
\end{proof}
\paragraph*{\bf{Remarks.}} \begin{itemize}\item It is easy to see that the above theorem
particularizes in the regular case ($G_{p_x}=G_x=0$) to the main
result of the reduced Energy-Momentum Method, (see page 35 in
\cite{SiLeMar}). In the regular case, one has to test the
stability of the Hessian $\ed_x^2 V_\mu$ restricted to
$(\g_\mu\cdot x)^\perp$, where $\mu$ is the momentum value of the
relative equilibrium under study and $V_\mu$ is Smale's
\emph{amended potential energy} (see \cite{Smale}), defined as
$$V_\mu (x)=V(x)+\frac 12
\langle\mu,\II^{-1}(x)(\mu)\rangle.$$ Obviously this function is
not well defined at those points $x\in Q$ such that $\dim \g_x
>0$, since the locked inertia tensor is not invertible. Hence, it
is not possible to define the Hessian of the amended potential in
the singular setting. However,  the following
relation is to be noted (see \cite{SiLeMar}, equation (2.28)), which holds in the
regular case at a relative equilibrium
$$\ed^2_x V_\mu (v_1,v_2) =\ed_x^2 V_\xi (v_1,v_2) + \langle
(\D\II\cdot v_1)(\xi),\II^{-1}(x)\left[(\D\II\cdot
v_2)(\xi)\right]\rangle .$$ This suggests that what we are testing
in Theorem \ref{REM} is exactly the singular analogue of the
Hessian of the amended potential, even when we cannot talk about
the amended potential itself.

\item In the regular case, according to our definition
\eqref{gsplitting} of $\mathfrak{t}$, $\Sigma=(\g_\mu\cdot
x)^\perp$ with respect to $\ll\cdot,\cdot\gg$. In the presence of
isotropy for the base point $x$ of our relative equilibria, it can
be seen that $\Sigma$ is the orthogonal complement to $\g_\mu\cdot
x$ within the space of admissible variations (see
\eqref{admissiblevariations} and the proof of
Proposition~\ref{rigidinternalisos}). Therefore, conditions $(ii)$
in Theorem \ref{REM} are tested in a space $\Sigma$ which is
orthogonally complementary (with respect to the kinetic energy
metric) to the drift orbit $G_\mu\cdot x$.
\end{itemize}
\section{An example: The sleeping Lagrange top}
We will apply our singular version of the reduced Energy-Momentum
Method to the study of the re\-la\-ti\-ve equilibrium known as the
\emph{sleeping Lagrange top}. This problem has been extensively
studied in the literature of Classical Mechanics; see in
particular \cite{LeRaSiMar} for a geometric perspective of the
problem using the Hamiltonian and symplectic formalism. Here we
show the advantages of the method stated in Theorem \ref{REM} when
studying stability in simple mechanical systems. Indeed, taking
into account the extra cotangent bundle structure of the problem
actually leads to simpler calculations for obtaining stability
results when compared with the non-adapted methods constructed for
general symmetric Hamiltonian systems.

 The Lagrange top is a symmetric simple mechanical
system defined on the configuration space $Q=\mathrm{SO}(3)$, in
the usual representation by orthogonal $3\times 3$ real matrices
with determinant 1, and being the symmetry group $G=\mathbb{T}^2$.
We use the right trivialization for $T\mathrm{SO}(3)$ given by the
isomorphism $T\mathrm{SO}(3)\rightarrow \mathfrak{so}(3)\times
\mathrm{SO}(3)$, i.e. $\delta g=\xi g$ with $\xi\in\mathfrak{g}$,
and identifying $\mathfrak{so}(3)$ with $\R^3$ under the inverse
of the usual isomorphism given by $\R^3\ni \mathbf{u}\rightarrow
\hat{\mathbf{u}}\in\mathfrak{so}(3)$ (as a $3\times 3$ matrix
algebra), defined by
$\hat{\mathbf{u}}\mathbf{v}=\mathbf{u}\times\mathbf{v},\,\forall\,
\mathbf{v}\in\R^3$. Analogous considerations hold to obtain the
trivialization of the phase space for the problem:
$T^*\mathrm{SO}(3)\rightarrow \mathbb{R}^3\times\mathrm{SO}(3)$.
Then we can write the Hamiltonian for this system as
\be\label{hamiltoniansleeping} H(\pi,\Lambda)=\frac 12 \pi\cdot
E_\Lambda^{-1}\pi+mgl\Lambda\mathbf{e_3}\cdot\mathbf{e_3},\ee
where $\Lambda\in Q,\, \pi\in\R^3$, and $(\pi,\Lambda)\in
T^*_\Lambda Q,\, E=\mathrm{diag}\, (i,i,i_3)$ and
$E_\Lambda=\Lambda E\Lambda^t$. Finally, $m,g,l$ are physical
constants of the model.

The group $\mathbb{T}^2$ is identified with $S^1\times S^1$ where,
in our matricial representation, both copies of the circle group
act by rotations around $\mathbf{e_3}$. Its action on $Q$ is given
by $(L,R)\cdot \Lambda=L\Lambda R^t$ and the action on the phase
space is by cotangent lifts. We do not actually need the (simple)
expression for the cotangent-lifted action since our methods rely
finally on the geometry of the action of $G$ on $Q$ once they were
constructed in this spirit. Let $(l,r)\in \R^2=
\mathrm{Lie}\,(\mathbb{T}^2)$, then the infinitesimal action of
$\R^2$ on $Q$ is given in the right trivialization of
$T\mathrm{SO}(3)$ by
$(l,r)_Q(\Lambda)=(l\mathbf{e_3}-r\Lambda\mathbf{e_3},\Lambda)$.

 To
finish with these preliminaries, let us show the Riemannian
structure of $Q$ associated to the kinetic energy of the problem.
If $V_1=(v_1,\Lambda)$ and $V_2=(v_2,\Lambda)$ are two tangent
vectors to $Q$ at $\Lambda$, then \be\label{metricsleeping}\ll
V_1,V_2\gg (\Lambda)=v_1\cdot E_\Lambda v_2\ee which is easily
checked to be a $\mathbb{T}^2$-invariant symmetric contravariant
bilinear tensor on $Q$. Note also from the
expression~\eqref{hamiltoniansleeping} of the Hamiltonian, that
the potential energy is given by
$V(\Lambda)=mgl\Lambda\mathbf{e_3}\cdot\mathbf{e_3}$.

It is well known that the phase space point $p_I=\FL((l,r)_Q(I))$,
where $I$ is the identity $3\times 3$ matrix, is a relative
equilibrium of this system, known as the sleeping Lagrange top.
Our aim is to study its stability. First of all, note that since
$\Lambda=I$, then $(l,r)_Q(I)=(\zeta\mathbf{e_3},I)$, where the
number $\zeta=l-r$ uniquely determines the infinitesimal
generator, and it is physically interpreted as the angular
velocity of the rigid body modelled by this system.

The first thing we need to compute is the locked inertia tensor.
It easily follows from the previous expressions that
$$\begin{array}{l} \ll (l_1,r_1)_Q(\Lambda),(l_2,r_2)_Q(\Lambda)\gg =
  (l_1\mathbf{e_3}-r_1\Lambda\mathbf{e_3})\cdot E_\Lambda
(l_2\mathbf{e_3}-r_2\Lambda\mathbf{e_3})\\$\quad$\\ \hspace{1.5cm}
= \left( l_1, r_1\right) \left(
\begin{array}{cc}\mathbf{e_3}\cdot E_\Lambda\mathbf{e_3} &
-\mathbf{e_3}\cdot\Lambda E\mathbf{e_3}\\
-\mathbf{e_3}\cdot E\Lambda^t\mathbf{e_3} & \mathbf{e_3}\cdot
E\mathbf{e_3}\end{array}\right)\left(
\begin{array}{c} l_2 \\ r_2
\end{array}\right)
\end{array}$$
and then $$\II(\Lambda)=\left(
\begin{array}{cc}\mathbf{e_3}\cdot E_\Lambda\mathbf{e_3} &
-\mathbf{e_3}\cdot\Lambda E\mathbf{e_3}\\
-\mathbf{e_3}\cdot E\Lambda^t\mathbf{e_3} & \mathbf{e_3}\cdot
E\mathbf{e_3}
\end{array}\right) .$$
Thus at the configuration $\Lambda=I$ of our relative equilibrium
we have
$$\II (I)=\left(
\begin{array}{cc}i_3 & -i_3 \\
-i_3 & i_3 \end{array}\right)$$ and so the momentum value of the
relative equilibrium is $$ \mu=\II (I)(l,r)=i_3(\zeta,-\zeta).$$
As the configuration point of our relative equilibrium is the
identity element of $\mathrm{SO}(3)$, then its isotropy group is
$H=G_I =S^1$, regarded as the diagonal embedding of $S^1$ in the
2-torus. For the momentum isotropy, just by noting that the group
is Abelian and thus the coadjoint representation is trivial, we
obtain $G_\mu=\mathbb{T}^2$. Finally, using the characterization
$G_{p_x}=G_x\cap G_\mu$ we also get $G_{p_x}=G_x=S^1$.

The next step is to choose a $G_{p_x}$-invariant complement to
$\g_{p_x}$ in $\g_\mu=\mathbb\R^2$ in order to implement the
splitting \eqref{optimal} and obtain the velocity $(l,r)^\perp$.
We will obtain it by orthogonality with respect to a
$G_{p_x}$-invariant inner product in $\g_\mu$. Since the adjoint
action is also trivial any inner product is $S^1$-invariant. We
can choose the family of inner products
$$\mathbf{G}=\left( \begin{array}{cc} k & 0\\ 0 &
1\end{array}\right).$$ Since $\g_{p_x}$ is generated by $(1,1)$,
its orthogonal complement $\mathfrak{p}$ with respect to
$\mathbf{G}$ is generated by the normalized vector
$\mathbf{k}=\frac{1}{\sqrt{k(1+k)}}(1,-k)$, so for the
infinitesimal generator of our relative equilibrium, $\xi=(l,r)$,
we have that $\xi^\perp$ is the orthogonal projection of $\xi$
onto $\mathfrak{p}$, i.e.
$\xi^\perp=\mathbf{G}(\xi,\mathbf{k})\mathbf{k}=\zeta
\left(\frac{1}{1+k},-\frac{k}{1+k}\right)$.\\
The augmented potential energy
$V_{\xi^\perp}(\Lambda)=V(\Lambda)-\frac 12 \II
(\Lambda)(\xi^\perp,\xi^\perp)$ corresponding to this
re\-la\-ti\-ve equilibrium is now written as $$
V_{\xi^\perp}(\Lambda)=mgl\Lambda\mathbf{e_3}\cdot\mathbf{e_3}-\frac{\zeta^2}{2(1+k)^2}\left[\mathbf{e_3}\cdot
E_\Lambda\mathbf{e_3}+2ki_3\mathbf{e_3}\cdot\Lambda\mathbf{e_3}+k^2i_3\right]
.$$ We now compute the derivative of the augmented potential.
 According to our right trivialization, a tangent vector $\delta
\Lambda\in T_\Lambda \mathrm{SO}(3)$ can be written as
$\delta\Lambda=\hat{\epsilon}\Lambda$, where $\epsilon$ is a
vector in $\R^3$. Then, after some vector calculus manipulations
we get
 \be \ed
V_{\xi^\perp}(\Lambda)\cdot\delta\Lambda =mgl\mathbf{e_3}\cdot
(\epsilon\times\Lambda\mathbf{e_3})-\frac{\zeta^2}{(1+k)^2}\left[ki_3\mathbf{e_3}\cdot(\epsilon\times\Lambda\mathbf{e_3})-
\mathbf{e_3}\cdot E_\Lambda(\epsilon\times\mathbf{e_3})\right] \ee
which vanishes at $\Lambda=I$, since $\FL ((l,r)_Q(I))$ is a
relative equilibrium. The Hessian of $V_{\xi^\perp}$ at
$\Lambda=I$ is
$$\begin{array}{l} \ed^2_{I}V_{\xi^\perp}(\delta
\Lambda_1,\delta\Lambda_2)=(\epsilon_1\times\mathbf{e_3})\cdot
(\epsilon_2\times\mathbf{e_3})
\left(\frac{(ki_3+i_3)\zeta^2}{(1+k)^2}-mgl\right)\\
\hspace{2cm}-\frac{\zeta^2}{(1+k)^2}(\epsilon_1\times\mathbf{e_3})\cdot
E(\epsilon_2\times\mathbf{e_3}) .\end{array}$$ A straightforward
computation shows that the correction term \eqref{correctionterm}
vanishes.

As the group is Abelian, then $\g_\mu=\g$ and $\q^\mu=(0,0)$, so
$\Sigma=\Sl$. Then by Theorem \ref{REM} the relative equilibrium
is $G_\mu$-stable if ${\ed^2_{I}V_{\xi^\perp}}\rrestr{\Sl}$ is
positive-definite. It is then necessary to obtain the linear slice
for the toral action on $\mathrm{SO}(3)$ at $\Lambda=I$ with
respect to the Riemannian metric~\eqref{metricsleeping}. A tangent
vector $\delta\Lambda\in T_{I}\mathrm{SO}(3)$ written as
$\delta\Lambda=\hat{\epsilon}I$ is orthogonal to $\g\cdot I$ if
and only if
$$\epsilon\cdot E\mathbf{e_3}=0$$
and thus we have $$\Sl=\{\delta\Lambda\in T_{I}\mathrm{SO}(3)\,
:\,\epsilon\in\mathrm{span}\{(1,0,0),(0,1,0)\}\}.$$ To check the
positive definiteness of ${\ed_I^2 V_{\xi^\perp}}\rrestr{\Sl}$ is
then equivalent to showing that the $2\times 2$ matrix
 \be\label{scalarmatrix}\left(\frac{(ki_3+i_3-i)\zeta^2}{(1+k)^2}-mgl\right)I_2
\ee is positive-definite, where $I_2$ is the identity matrix in
$\R^2$; so the condition for stability is satisfied if
\be\label{ltopcond}\zeta^2>\frac{(1+k)^2mgl}{ki_3+i_3-i}\ee
\paragraph{\bf{Remark.}}The above expression
is exactly the one obtained in \cite{OrRa99}, page 718, by using
the (singular) Energy-Momentum Method for general Hamiltonian
systems in arbitrary symplectic manifolds, i.e.
Theorem~\ref{lermansinger} in this paper. In that work, the same
condition for the $\mathbb{T}^2$-stability of the sleeping
Lagrange top is obtained after computing algebraically the
eigenvalues of a $4\times 4$ matrix
($\ed_{p_x}^2h_\xi\rrestr{V_s}$), which was not put in
 block-diagonal form by applying their general stability criterion.
 With our method the same expression follows easily
from the unique eigenvalue of the scalar matrix
\eqref{scalarmatrix}. This shows the potential power of employing
methods adapted to the cotangent bundle structure in the study of
simple mechanical systems.  In particular in practical stability
problems with higher dimensional configuration spaces the
implementation of Theorem~\ref{lermansinger} could involve
checking the definiteness of much more complicated matrices,
for\-cing the use of numerical methods in some situations where
the reduced Energy-Momentum Method developed in this paper
(Theorem \ref{REM} and also Corollary \ref{blockdiagonalREM})
could offer simpler or even exact results.

For the sake of completeness we will sharpen the stability
condition \eqref{ltopcond}, follo\-wing \cite{OrRa99}. Since in
the stability condition $k$ appears, which is related to
$\mathbf{G}$, the sharpest (or optimal) stability condition will
be the lowest value of $\zeta^2$ among all possible $k$. This is
an straightforward optimization problem in elementary calculus, so
by differentiating the expression
$$f(k)=\frac{(1+k)^2mgl}{ki_3+i_3-i} ,$$
we find that the minimum value is reached at
$k=\frac{2i-i_3}{i_3}$. So the sharpest stability condition yields
the well-known lower bound for the angular velocity
$$\zeta^2>\frac{4mgli}{i_3^2}.$$
\section{The singular Arnold form and Block-Dia\-go\-na\-li\-za\-tion}
\label{singulararnoldsection} In the previous section we used the
realization of the symplectic normal space $V_s$ at a relative
equilibrium given in \eqref{Vre} and which has been shown to be
isomorphic to $\q^\mu\oplus\Sl\oplus\Sl^*$ by the isomorphism
$\eqref{isore2}$. This was helpful to simplify the study of the
definiteness of ${\ed_{p_x}^2h_{\xi^\perp}}\rrestr{V_s}$,
obtaining a block-diagonal structure with one block being
positive-definite or trivial, reducing the problem to study the
definiteness of $\left(\ed^2_xV_{\xi^\perp}+\corr_{\xi^\perp}
(x)\right)\rrestr{\Sigma}$. In this section we pursue the study of
the symplectic normal space $V_s$ by obtaining new block-diagonal
forms for ${\ed_{p_x}^2h_{\xi^\perp}}\rrestr{V_s}$ and also for
the symplectic matrix $\Omega$ at a relative equilibrium. Recall
that in the $I$-representation, the symplectic matrix $\w (p_x)$
of $T_{p_x}(T^*Q)$ has the following form (see \cite{Paternain}
and \cite{PerRoSD2}). If
$(\eta_1,a_1;\,\nu_1,\alpha_1),\,(\eta_2,a_2;\,\nu_2,\alpha_2)\in
I(T_{p_x}(T^*Q))=(\lir\oplus\Sl)\oplus (\lir^*\oplus\Sl^*)$ then
\be\label{Isymplecticmatrix} \Omega
\left((\eta_1,a_1;\,\nu_1,\alpha_1),(\eta_2,a_2;\,\nu_2,\alpha_2)\right)=
\langle \nu_2,\eta_1\rangle+\langle\alpha_2,a_1\rangle-\langle
\nu_1,\eta_2\rangle-\langle\alpha_1,a_2\rangle .\ee
 In the following we will use a parametrization of $V_s$ different from \eqref{isore2}. This, when available,
 will present the extra advantage of putting
${\ed_{p_x}^2h_{\xi^\perp}}\rrestr{V_s}$ in block-diagonal form
consisting of three blocks instead of two, as in
Theorem~\ref{REM}, hence simplifying the stability
analysis further.
 We will start by
defining a singular analogue of the Arnold form, (see
\cite{Arnold2} and \cite{MarLec} to see how the Arnold form arises
in the study of the stability of regular relative equilibria of
simple mechanical systems). Hereafter, we assume that we will be working
under the same conditions and hypotheses as in the previous
section. In particular, we will always have $p_x=\FL(\xi_Q(x))$ as
a relative equilibrium of the simple mechanical system
\eqref{smshamiltonian} with momentum $\mu$ and base point isotropy
$G_x=H$.
\begin{defn}\label{arnold}
The singular Arnold form at a relative equilibrium with base point
$x$ and momentum $\mu$ is the bilinear form on $\q^\mu$, $\Ar
:\q^\mu\times\q^\mu\rightarrow \R$ defined by $$ \Ar
(\lambda_1,\lambda_2)=\langle \ad^*_{\lambda_1}\mu,\Lambda
(x,\mu)(\lambda_2)\rangle ,$$ where the map $\Lambda
(x,\mu):\q^\mu\rightarrow \lir$ is defined by $$\Lambda
(x,\mu)(\lambda)=\rII^{-1}\left(
\ad^*_\lambda\mu\right)+\Proj_{\lir^*}\left[\ad_\lambda
\left(\rII^{-1}\mu\right)\right] .$$
\end{defn}
\paragraph{\bf{Remark.}} To the best of our knowledge, the first time a singular
analogue of the Arnold form appears in the literature is in
\cite{Lewisblock}, equation (3.56), under the name of
``generalized Arnold form''. In that work, this object is defined
in the context of general Lagrangian systems and the Lagrangian
Block-Diagonalization method. Also see Section \ref{remarks} for a
comparison of other results in  \cite{Lewisblock}.

We will now define another two spaces which will be useful for the
obtention of block-diagonal expressions. The motivation for the
introduction of these spaces is, following the ideas in
\cite{SiLeMar}, that provided a non-degeneracy condition is
satisfied, ${\ed_{p_x}^2h_{\xi^\perp}}\rrestr{V_s}$ will further
block-diagonalized, refining the conditions of Theorem~\ref{REM}.

 Recalling the
identification $\kappa:\q^\mu\oplus \Sl\oplus \Sl^*\rightarrow
V_s$, let us define the following subspace of $\q^\mu\oplus \Sl$:
\be\label{wintchar} w_\text{int} = \left\{(\lambda^b,b)\in
\q^\mu\oplus \Sl\, :\, (\D\II\cdot
(\lambda^b_Q(x)+b))(\xi^\perp)\in\mathfrak{p}^*\right\} .\ee Also,
using the map $(\lambda,b)\mapsto \lambda_Q(x)+b$ which maps
isomorphically $\q^\mu\oplus\Sl$ onto $\Sigma\subset T_xQ$, we
define the following subspaces of $\Sigma$: \be\label{sigmarig}
\Sigma_\text{rig}
 =  \left\{\lambda_Q(x)\, :\,\lambda\in\q^\mu\right\}
,\quad\mathrm{and}\ee \be\label{sigmaint} \Sigma_\text{int}  =
\left\{\lambda^b_Q(x)+b\, :\, (\lambda^b,b)\in
w_\text{int}\right\}. \ee Note that we have the following obvious
identifications:
$$\begin{array}{lll}  \Sigma_\text{rig}  & \simeq &
\q^\mu \\
 \Sigma_\text{int} & \simeq & w_\mathrm{int}.
\end{array}$$
Following \cite{SiLeMar} we can give the following interpretation
for these two spaces: $\Sigma$ is seen as the space of all
admissible variations orthogonal to the infinitesimal drift
directions $\g_\mu\cdot x$. This space has a contribution
$\Sigma_\text{rig}$ which corresponds to variations in $\Sigma$
which generate group motions, i.e. regarding our systems as a
``rigid body" without internal structure. On the contrary, the
subspace $\Sigma_\text{int}$ corresponds to all the variations of
our system that are purely internal, i.e. variations in ``shape",
not coming from the symmetry group.
\begin{prop}\label{rigidinternalisos}
If the Arnold form is non-degenerate then $\Sigma =
\Sigma_\mathrm{rig}\oplus\Sigma_\mathrm{int}.$
\end{prop}
For the proof of this proposition, we will need the following
lemma
\begin{lemma}\label{admissible}
For every $\lambda\in\q^\mu$ and $v\in\Sl$\begin{itemize}
\item[(i)]
$\Proj_\h\,\left[(\D\II\cdot\lambda_Q(x))(\xi^\perp)\right]=0$
\item[(ii)] $\Proj_\h\,\left[(\D\II\cdot v)(\xi^\perp)\right]=0$.
\end{itemize}
\end{lemma}
\begin{proof}
For (i), using \eqref{propertieslocked2}, for any $\zeta\in\h$ one
has
$$\begin{array}{lll}
(\D\II\cdot\lambda_Q(x))(\xi^\perp,\zeta) & = & -\II
(x)(\ad_\lambda\xi^\perp,\zeta)-\II(\xi^\perp,\ad_\lambda\zeta)\\
& = & -\langle\ad^*_\lambda\mu,\zeta\rangle =  0 ,
\end{array}$$
where the second equality follows since $\ker\II (x)=\h$ and
$\II(x)(\xi^\perp)=\mu$.\\
 For (ii), making $\lambda=\zeta\in\h$ in
item (iii) of Theorem \ref{covariant}, we have $(\D\II\cdot
v)(\xi^\perp,\zeta)=0$ for every $\zeta\in\h$.
\end{proof}
\begin{proof}(of the proposition) It is clear that
$\Sigma_\mathrm{rig}\subset \Sigma$ and
$\Sigma_\mathrm{int}\subset\Sigma$ so we have to prove
$\Sigma_\mathrm{rig}\cap\Sigma_\mathrm{int}=0$ and
$\Sigma_\mathrm{rig}+\Sigma_\mathrm{int}=\Sigma$. Let
$0\neq\lambda\in\q^\mu$. Then $\lambda_Q(x)\in
\Sigma_\mathrm{rig}\cap\Sigma_\mathrm{int}$ if and only if
$(\D\II(x)\cdot\lambda_Q(x))(\xi^\perp,\epsilon)=0$ for every
$\epsilon \in \mathfrak{t}+\h$. By $(i)$ in Lemma
\ref{admissible}, this holds if and only if the same condition is
satisfied for every $\epsilon\in\mathfrak{t}$. Using
\eqref{propertieslocked2}, this is equivalent to
$$\begin{array}{lll}0 & = & \II (x)(\ad_\lambda\xi^\perp,\epsilon)+\II
(x)(\xi^\perp,\ad_\lambda\epsilon) =  \langle
\II(x)\left(\ad_\lambda\left(
\rII^{-1}\mu\right)\right)+\ad^*_\lambda\mu,\epsilon\rangle\\
& = & \langle \rII\left(\Proj_{\lir^*}\left[\ad_\lambda\left(
\rII^{-1}\mu\right)\right]\right)+\ad^*_\lambda\mu,\epsilon\rangle
= \langle \rII\left(\Proj_{\lir^*}\left[\left(\ad_\lambda
\rII^{-1}\mu\right)\right]+\rII^{-1}(\ad^*_\lambda\mu
)\right),\epsilon\rangle
\end{array}$$
for every $\epsilon\in\mathfrak{t}$, regarding that
$\ad_\lambda^*\mu\in\lir$ if $\lambda\in\q^\mu$. Since $\rII$ is
an isomorphism, this condition is the same as $\Lambda
(x,\mu)(\lambda)\in\mathfrak{p}\subset\g_\mu$, but then the Arnold
form would be degenerate, which is a contradiction.

To prove that $\Sigma_\mathrm{rig}+ \Sigma_\mathrm{int}=\Sigma$
let us note the following: if we call
$$\mathcal{D}=\left\{\delta q\in T_x Q\, :\,
\Proj_\h\,\left[(\D\II\cdot \delta
q)(\xi^\perp)\right]=0\right\},$$ then by the definitions of
$\mathfrak{t}$ \eqref{gsplitting}, $\q^\mu$ \eqref{qmu}, and by
Lemma \ref{admissible} we have that $\Sigma$ is precisely the
orthogonal complement to $\g_\mu\cdot x$ in $\mathcal{D}$. The
rest of the proof is then a consequence of Proposition 3.7 in
\cite{Lewisblock}.
\end{proof}
\begin{cor}\label{isoarnoldcor}
If the Arnold form is non-degenerate then the map
$\tilde{\kappa}:\q^\mu \oplus \Sigma_\mathrm{int}\oplus
\Sl^*\rightarrow V_s$ defined as
\be\label{ktilde}\tilde{\kappa}(\lambda,(\lambda^a_Q(x)+a),\gamma)=\kappa
(\lambda+\lambda^a,a,\gamma)\ee for every $\lambda\in \q^\mu,\,
(\lambda^a,a)\in w_\mathrm{int}$ and $\gamma\in \Sl^*$ is a
$G_{p_x}$-equivariant isomorphism.
\end{cor}
 We will assume now that the Arnold form is non-degenerate and
 then we
will study the symplectic matrix $\Omega$ of $V_s$ and
$\ed^2_{p_x}h_{\xi^\perp}$.
\begin{prop}[Block-Diagonalization
forms]\label{blockdiagonalforms} If the Arnold form is
non-degenerate, under the isomorphism $\tilde{\kappa}$ of
Corollary \ref{isoarnoldcor} we have the following expressions for
the symplectic matrix $\Omega$ and
${\ed_{p_x}^2h_{\xi^\perp}}\rrestr{V_s}$:
\be\label{blocksigma}\begin{array}{ccccc}
 & & \q^\mu & \Sigma_\mathrm{int} & \Sl^*\\
\Omega & = & \left(\begin{array}{c} \Xi \\ \Psi^t \\ 0 \end{array}\right. & \begin{array}{c} -\Psi \\
S_\mu \\ -\mathbf{1} \end{array} & \left. \begin{array}{c} 0 \\
\mathbf{1}
\\ 0 \end{array}\right)
\end{array}\ee
and \be\label{blockhess}\begin{array}{ccccc}
 & & \q^\mu & \Sigma_\mathrm{int} & \Sl^*\\
{\ed_{p_x}^2h_{\xi^\perp}}\rrestr{V_s} & = &
\left(\begin{array}{c} \mathrm{Ar}
\\ 0 \\ 0 \end{array}\right.
& \begin{array}{c} 0 \\
 (\ed^2_{x}V_{\xi^\perp}+\corr_{\xi^\perp}
(x))\rrestr{\Sigma_\mathrm{int}}\\ 0 \end{array} & \left. \begin{array}{c} 0 \\
0
\\ \ll\cdot,\cdot\gg_{\Sl^*} \end{array}\right)
\end{array},\ee
where the entries of $\Omega$ are $$\begin{array}{rll} \Xi
(\lambda_1,\lambda_2) & = & -\langle
\mu,\ad_{\lambda_1}\lambda_2\rangle\vspace{1mm}\\ \Psi
(\lambda,(\lambda^b_Q(x)+b)) & = & \langle\mu,\ad_\lambda
\lambda^b\rangle\vspace{1mm}\\
S_\mu ((\lambda^a_Q(x)+a),(\lambda^b_Q(x)+b)) & = & -\langle
\mu,\ad_{\lambda^a}\lambda^b\rangle+\ll C(b)(\xi^\perp),a\gg_\Sl\\
& & - \ll C(a)(\xi^\perp),b\gg_\Sl .\end{array}$$
\end{prop}
\begin{proof}
The form for $\Omega$ follows trivially from
\eqref{Isymplecticmatrix} and the definition of $\tilde{k}$ from
\eqref{ktilde}.
 For ${\ed_{p_x}^2h_{\xi^\perp}}\rrestr{V_s}$, and recalling its block-diagonal form
showed in the proof of Theorem~\ref{REM}, the only two things that
we must check are
\begin{itemize}
\item[(i)] $(\ed^2_xV_{\xi^\perp}+\corr_{\xi^\perp}
(x))({\lambda_1}_Q(x),{\lambda_2}_Q(x))=\Ar
(\lambda_1,\lambda_2)$, and \item[(ii)]
$(\ed^2_xV_{\xi^\perp}+\corr_{\xi^\perp} (x))(\lambda_Q(x),\delta
q)=0$,
\end{itemize}
for  $\lambda,\lambda_1,\lambda_2\in \q^\mu$ and $\delta q\in
\Sigma_\mathrm{int}$.

 To prove (i) recall that the potential
energy $V$ is $G$-invariant and then
$$(\ed^2_xV_{\xi^\perp}+\corr_{\xi^\perp} (x))({\lambda_1}_Q(x),{\lambda_2}_Q(x))=
\left( -\frac 12\ed_x^2 \left(\II(x) (\xi^\perp,\xi^\perp)\right)
+\corr_{\xi^\perp} (x)
\right)({\lambda_1}_Q(x),{\lambda_2}_Q(x)).$$ We compute now both
terms in the right hand side of the above expression. For the
first one,
$$\begin{array}{lll}
 -\frac 12 \ed_x^2 \left(\II(x) (\xi^\perp,\xi^\perp)\right)({\lambda_1}_Q(x),{\lambda_2}_Q(x)) & = &
 -\frac 12  {\lambda_2}_Q((\D\II\cdot {\lambda_1}_Q(x))(\xi^\perp,\xi^\perp))
 (x)\vspace{0.5mm}\\
 & = & {\lambda_2}_Q(\II
 (\ad_{\lambda_1}\xi^\perp,\xi^\perp))(x)\vspace{0.5mm}\\
 & = & (\D\II\cdot
 {\lambda_2}_Q(x))(\ad_{\lambda_1}\xi^\perp,\xi^\perp)\vspace{0.5mm}\\
 & = & -\II (x) (\ad_{\lambda_2}(\ad_{\lambda_1}\xi^\perp),\xi^\perp)\vspace{0.5mm}\\
 &  & -\II (x)(\ad_{\lambda_1}\xi^\perp,\ad_{\lambda_2}\xi^\perp).
\end{array}$$\\
 Now for the second,
$$\begin{array}{lll}
\corr_{\xi^\perp} (x)({\lambda_1}_Q(x),{\lambda_2}_Q(x)) & = &
\langle\Proj_\lir\left[(\D\II\cdot
{\lambda_1}_Q(x))(\xi^\perp)\right],
\rII^{-1}\Proj_\lir\left[(\D\II\cdot {\lambda_2}_Q(x))(\xi^\perp)\right]\rangle\vspace{1mm} \\
& = & \langle \II (x)(\ad_{\lambda_1}\xi^\perp)+
\ad^*_{\lambda_1}\mu,\rII^{-1}\left[\II (x)(\ad_{\lambda_2}\xi^\perp)+\ad^*_{\lambda_2}\mu  \right]\rangle\vspace{1mm}\\
& = & \II (x)(\ad_{\lambda_2}\xi^\perp,\ad_{\lambda_1}\xi^\perp)+\langle\ad^*_{\lambda_1}\mu,\ad_{\lambda_2}\xi^\perp\rangle\vspace{1mm}\\
 &  & +
\langle\ad^*_{\lambda_2}\mu,\ad_{\lambda_1}\xi^\perp\rangle
+\langle\ad^*_{\lambda_1}\mu,\rII^{-1}(\ad^*_{\lambda_2}\mu
)\rangle,
\end{array}$$
where we have used that $\II (x)(\eta_1,\eta_2)=\II
(x)(\eta_1^\lir,\eta_2^\lir)$ $\forall\, \eta_1,\eta_2\in\g ,$
since $\ker \II (x)=\h$. Also,
$\langle\ad_\lambda^*\mu,\eta\rangle=\langle\ad_\lambda^*\mu,\eta^\lir\rangle$
$\forall\,\eta\in\g ,$ since $\lambda\in\q^\mu$, which means that
$\ad_\lambda^*\mu\in\h^\circ$. Finally, noticing that
$$\langle\ad^*_{\lambda_2}\mu,\ad_{\lambda_1}\xi^\perp\rangle =\II (x)(\xi^\perp,\ad_{\lambda_2} (\ad_{\lambda_1}\xi^\perp)),$$
and putting both contributions together, we obtain
$$\begin{array}{lll}
(\ed^2_xV_{\xi^\perp}+\corr_{\xi^\perp}
(x))({\lambda_1}_Q(x),{\lambda_2}_Q(x)) & = &
\langle\ad_{\lambda_1}^*\mu,\ad_{\lambda_2}\xi^\perp
+\rII^{-1}(\ad^*_{\lambda_2}\mu)\rangle\vspace{0.6mm}\\
& = &
\langle\ad_{\lambda_1}^*\mu,\Proj_\lir\,\left[\ad_{\lambda_2}\left(\rII^{-1}\mu\right)\right]
+\rII^{-1}(\ad^*_{\lambda_2}\mu)\rangle\\
 & = & \Ar (\lambda_1,\lambda_2),
\end{array}$$
since $\ad_{\lambda_1}^*\mu\in\h^\circ$ and
$\xi^\perp=\rII^{-1}\mu$.

Again for (ii) because the potential energy $V$ is $G$-invariant
we have
$$\begin{array}{l} \left(\ed^2_xV_{\xi^\perp}+\corr_{\xi^\perp}
(x)\right)(\lambda_Q(x) ,\lambda^a_Q(x) +a)=\vspace{0.5mm}\\ =
\left( -\frac 12\ed_x^2 \left(\II(x) (\xi^\perp,\xi^\perp)\right)
+\corr_{\xi^\perp} (x)
\right)(\lambda_Q(x),\lambda^a_Q(x)+a).\end{array}$$ We will start
by computing the contribution of the correction term:
$$\begin{array}{lll}
\corr_{\xi^\perp} (x)(\lambda ,\lambda^a_Q(x)+a) & = & \langle
\Proj_\lir\left[(\D\II\cdot\lambda_Q(x))(\xi^\perp)
\right],\rII^{-1}\Proj_\lir\left[(\D\II\cdot
(\lambda^a_Q(x)+a))(\xi^\perp) \right]\rangle\vspace{1mm} \\
& = & -\langle\Proj_\lir\left[\II (x)(\ad_\lambda\xi^\perp)
\right]+\ad^*_\lambda\mu,\rII^{-1}\Proj_\lir\left[(\D\II\cdot
(\lambda^a_Q(x)+a))(\xi^\perp) \right]\rangle\vspace{1mm} \\
& = & -\langle\II
(x)(\ad_\lambda\xi^\perp)+\ad^*_\lambda\mu,\rII^{-1}(\D\II\cdot
(\lambda^a_Q(x)+a)(\xi^\perp)\rangle ,
\end{array}$$
since $\ker \II (x)=\h$ and by Lemma \ref{admissible}
$\Proj_\h\,\left[(\D\II\cdot
(\lambda^a_Q(x)+a)(\xi^\perp)\right]=0.$ For the first term we
have
$$-\frac
12\ed_x^2 \left(\II(x) (\xi^\perp,\xi^\perp)\right)(\lambda
,\lambda^a_Q(x)+a)  =  (\D\II\cdot
(\lambda^a_Q(x)+a))(\ad_\lambda\xi^\perp,\xi^\perp),$$ and so we
finally obtain
$$\left(\ed^2_xV_{\xi^\perp}+\corr_{\xi^\perp}
(x)\right)(\lambda_Q(x) ,\lambda^a_Q(x)
+a)=-\langle\ad^*_\lambda\mu,\rII^{-1}(\D\II\cdot
(\lambda^a_Q(x)+a)(\xi^\perp)\rangle .$$ This expression is zero
since by construction $(\D\II\cdot (\lambda^a_Q(x)+a))(\xi^\perp)$
annihilates $\mathfrak{t}$ if $(\lambda^a,a)\in w_\mathrm{int}$
and $\rII^{-1}(\ad^*_\lambda\mu)\in\mathfrak{t}$ for every
$\lambda\in\q^\mu$. To see this, note that the image of
$\rII^{-1}$ is in $\lir$ and hence $\II
(x)(\rII^{-1}(\ad^*_\lambda\mu),\zeta)=\langle\ad^*_\lambda\mu,\zeta\rangle=-\langle\ad^*_\zeta\mu,\lambda\rangle=0,$
for every $\zeta\in\g_\mu$, in particular if
$\zeta\in\mathfrak{p}$.
\end{proof}
An inspection of the form of
${\ed_{p_x}^2h_{\xi^\perp}}\rrestr{V_s}$ in \eqref{blockhess}
together with Theorem \ref{REM} leads to the following sharper
result concerning the $G_\mu$-stability of $p_x$.
\begin{cor}[Block-diagonalization and
stability]\label{blockdiagonalREM} In the hypothesis of Theorem
\ref{REM}, and assuming that the Arnold form in non-degenerate,
if:
\begin{itemize} \item[(i)] $\dim Q-\dim
G+\dim G_x >0$ \item[(ii)]
$\left(\ed^2_{x}V_{\xi^\perp}+\corr_{\xi^\perp}
(x)\right)\rrestr{\Sigma_\mathrm{int}}$ is positive definite
\item[(iii)] the Arnold form is positive definite,
\end{itemize}
or
\begin{itemize} \item[(i)] $\dim Q-\dim
G+\dim G_x =0$ \item[(ii)] the Arnold form is definite (positive
or negative),
\end{itemize}
then the relative equilibrium is $G_\mu$-stable.
\end{cor}
\paragraph{\bf{Remark.}} In \cite{SiLeMar}, Theorem 2.7, equivalent results
 to our Proposition~\ref{blockdiagonalforms} and
Corollary~\ref{blockdiagonalREM} are obtained as a consequence of
their reduced Energy-Momentum Method, for the particular case
$G_x=\{e\}$ (regular relative equilibria).
\section{Some remarks on the stability results}\label{remarks}
\paragraph{\bf{The residual symmetry sub-blocking.}}$\quad$\\
In this subsection we will use the fact that the symplectic normal
space $V_s$ supports a linear representation of $G_{p_x}$ in order
to improve our block-diagonalization results. We start with the
following lemma (see \cite{PerRoSD2} for a proof):
\begin{lemma}\label{intermch3}
Let $h\in H=G_x,\,\xi\in\lir$ and $v,w\in T_xQ$, then
\begin{itemize}\item[(i)] $ \Ad_{h^{-1}}^*\left[(\D\II\cdot v)(\xi)\right] =
(\D\II\cdot (h\cdot v))(\Ad_h\xi)$ \item[(ii)] $\ll C(h\cdot
v)(\Ad_{h}\xi),h\cdot w\gg_\Sl=\ll C(v)(\xi),w\gg_\Sl$.
\end{itemize}
\end{lemma}
As an immediate consequence of this and just by regarding their
definitions \eqref{sigmadef}, \eqref{qmu}, \eqref{sigmaint}, the
spaces $\Sigma, \q^\mu$ and $\Sigma_\text{int}$ are
$G_{p_x}$-invariant. We will use a tool from representation theory
known as the \emph{isotypic decomposition} of a linear space acted
linearly upon a compact Lie group to take advantage of the
residual symmetry group $G_{p_x}$ in order to further
block-diagonalize ${\ed_{p_x}^2h_{\xi^\perp}}\rrestr{V_s}$. For
this, we need the following definition, which is taken from
\cite{Golu}.
\begin{defn}
Let $K$ be a compact Lie group acting li\-nearly on a (real and
finite dimensional) linear space $N$. The isotypic decomposition
of $N$ is the unique decomposition
$$N=N_1\oplus\ldots\oplus N_r,$$
where each ${N}_i$ is the direct sum of all $K$-isomorphic
irreducible subspaces of $N$, and it is called an \emph{isotypic
component} of $N$.
\end{defn}
The  isotypic decomposition of a linear space satisfies the
following remarkable property (see \cite{Golu}): If $B$ is a
$K$-invariant bilinear form on $N$ (that is, $B(g\cdot v_1,g\cdot
v_2)=B(v_1,v_2)$ for every $g\in K,\,v_1,v_2\in N$) then
$B\rrestr{{N}_i\times {N}_j}=0$ for every pair of isotypic
components of $N$ with ${N}_i\neq {N}_j$. Therefore the expression
of $B$ block-diagonalizes with respect to the isotypic
decomposition of $N$. We will apply this property to the bilinear
form given by $\ed^2_xV_{\xi^\perp}+\corr_{\xi^\perp}(x)$.
\begin{lemma}
Let $p_x$ be a relative equilibrium of the simple mechanical
system \eqref{smshamiltonian} with velocity $\xi$. Then the
bilinear form $\ed^2_xV_{\xi^\perp}+\corr_{\xi^\perp}(x)$ is
$G_{p_x}$-invariant.
\end{lemma}
\begin{proof}We have to prove separately the invariance of each
term. For $\ed^2_xV_{\xi^\perp}$ the result follows if we prove
that $V_{\xi^\perp}$ is $G_{p_x}$-invariant. Since $V$ is
$G$-invariant, and so $G_{p_x}$-invariant, we only need to prove
that the function $\frac 12\II (\cdot)(\xi^\perp,\xi^\perp)$ is
$G_{p_x}$-invariant. Recall that at a relative equilibrium,
besides the characterization $G_{p_x}=G_\mu\cap G_x$, one also has
$G_{p_x}=\{h\in G_x\, :\, \Ad_h\xi^\perp=\xi^\perp\},$ which
follows trivially from the property
$$g\in G_{p_x}\Leftrightarrow g\cdot\FL (\xi^\perp_Q(x))=\FL (\xi^\perp_Q(x)),$$
since $\xi^\perp\in\mathfrak{p}\subset\lir$ and by definition our
relative equilibrium can be written as $p_x=\FL
(\xi_Q(x))=\FL(\xi^\perp_Q(x))$. Recall also the invariance
property of the locked inertia tensor \eqref{propertieslocked1}.
Then, for every $x'\in Q$ and $h\in G_{p_x}$, one has
$$\frac 12\II (h\cdot x')(\xi^\perp,\xi^\perp)=\frac 12\II
(x')(\Ad_{h^{-1}}\xi^\perp,\Ad_{h^{-1}}\xi^\perp)=\frac 12\II
(x')(\xi^\perp,\xi^\perp).$$
 For the correction term, recall from
Lemma~\ref{intermch3}, that if $h\in G_{p_x}$ and $\delta q\in
T_xQ$, then $$(\D\II\cdot (h\cdot\delta
q))(\xi^\perp)=\Ad^*_{h^{-1}}\left( (\D\II\cdot \delta
q)(\xi^\perp) \right),$$ since $\Ad_h\xi^\perp=\xi^\perp$. Note
also that $\rII:\lir\rightarrow \lir^*$ is a $G_x$-equivariant
isomorphism, that is $\rII\circ \Ad_{h}=\Ad_{h^{-1}}^*\circ \rII
.$ Then, given $\delta q_1,\delta q_2\in T_xQ$ and $h\in G_{p_x}$,
we have
$$\begin{array}{ll}
\corr_{\xi^\perp}(x)(h\cdot\delta q_1,h\cdot\delta q_2) & =
\langle(\D\II\cdot (h\cdot\delta
q_1))(\xi^\perp),\rII^{-1}\left((\D\II\cdot
(h\cdot\delta q_2))(\xi^\perp) \right)\rangle\\
& = \langle\Ad^*_{h^{-1}}\left( (\D\II\cdot \delta q_1)(\xi^\perp)
\right),\rII^{-1}\left( \Ad^*_{h^{-1}}\left( (\D\II\cdot \delta
q_2)(\xi^\perp)\right)\right)\rangle\\
& = \langle\Ad^*_{h^{-1}}\left( (\D\II\cdot \delta q_1)(\xi^\perp)
\right),\Ad_{h}\left( \rII^{-1}\left( (\D\II\cdot \delta
q_2)(\xi^\perp)\right) \right)\rangle\\
& = \langle (\D\II\cdot\delta q_1)(\xi^\perp),\rII^{-1}\left(
(\D\II\cdot \delta q_2)(\xi^\perp) \right)\rangle\\
& = \corr_{\xi^\perp}(x)(\delta q_1,\delta q_2)
\end{array}$$\end{proof}
Thus, in the hypothesis of Theorem~\ref{REM} we can use the fact
that $\Sigma$ is a $G_{p_x}$-invariant subspace of $T_xQ$, and
then testing the definiteness of
$(\ed^2_xV_{\xi^\perp}+\corr_{\xi^\perp}(x))\rrestr{\Sigma}$ is
equivalent to testing the definiteness of every restriction
$(\ed^2_xV_{\xi^\perp}+\corr_{\xi^\perp}(x))\rrestr{\Sigma_i}$
where
$$\Sigma=\Sigma_1\oplus\ldots\oplus \Sigma_r$$
is the isotypic decomposition of $\Sigma$. Analogously, if the
Arnold form is not degenerate, to test the conditions for
stability of Corollary~\ref{blockdiagonalREM} is equivalent to
test definiteness of $\Ar_{\q^\mu_i}$ and
$(\ed^2_xV_{\xi^\perp}+\corr_{\xi^\perp}(x))\rrestr{{\Sigma_\mathrm{int}}_i}$
for each of the isotypic components of $\q^\mu$ and
$\Sigma_\mathrm{int}$, respectively.
\paragraph{\bf{Nature of the stability results.}}$\quad$\\
In Theorems \ref{lermansinger} and \ref{REM} we have imposed the
compactness of the momentum isotropy subgroup $G_\mu$. However,
this compactness condition can be weakened by assuming that $\mu$
is \emph{split}, and this is how the original theorems are stated
in \cite{LeSi} and \cite{OrRa99}. A momentum value $\mu$ is called
split if there exists a $G_\mu$-invariant complement to $\g_\mu$
in $\g$. Obviously this is the case if $G_\mu$ is compact, since
in this case one can define this complement to be the orthogonal
complement to $\g_\mu$ with respect to some $G_\mu$-invariant
inner product on $\g$. Likewise, if $G$ itself is compact or
Abelian, then every momentum value is automatically split.

In the most general situation, if the relative equilibrium under
study has not a split momentum value, Theorem~\ref{REM} and
Corollary~\ref{blockdiagonalREM} are still applicable, but in that
case one does not obtain conditions for $G_\mu$-stability, only
for the weaker notion of \emph{leafwise stability}. A relative
equilibrium is called leafwise stable if it is $G_\mu$-stable for
the restriction of the Hamiltonian flow to $\J^{-1}(\mu)$. The
reason for this nomenclature is that in the free case, a relative
equilibrium $z$ with momentum $\mu$ for the symmetric Hamiltonian
system $(\PP,\w,G,\J, h)$ is leafwise stable if the point $[z]\in
\PP/G$ in the orbit space is a Lyapunov stable equilibrium for the
reduced Hamiltonian system on the symplectic leave
$\J^{-1}(\mu)/G_\mu\subset \PP/G$, rather than being stable in the
full Poisson quotient $\PP/G$. The results of 
\cite{LeSi,OrRa99,Patrick} guarantee that if $\mu$ is split,
then leafwise stability of $z$ implies Lyapunov stability of $[z]$
in $\PP/G$, and, hence, so do Theorem~\ref{REM} and
Corollary~\ref{blockdiagonalREM}. See \cite{PatRobWul} for a more
detailed explanation of these concepts.
\paragraph{\bf{The reduced Energy-Momentum Method and Lagrangian Block-\\ Diagonalization.}}$\quad$\\
In \cite{Lewisblock}  a method was constructed for testing the
stability of relative equilibria of symmetric Lagrangian systems.
In that work, the techniques of the reduced Energy-Momentum Method
of \cite{SiLeMar} are translated to systems defined on the tangent
bundle $TQ$ of the configuration space and developed for general
Lagrangian systems invariant under a possibly non free, tangent
lifted action. We briefly explain the relationship of the results
of \cite{Lewisblock} applied to simple mechanical systems with our
work. See \cite{Lewisblock,Lewisdrop} for more details
on the Lagrangian Block-Diagonalization method.

Let $L\in C^G (TQ)$ be a function on the tangent bundle of $Q$
invariant under the tangent lift of a proper action of the Lie
group $G$ on $Q$. This function is called a Lagrangian. There is a
well-known procedure to obtain a bundle map $\FL:TQ\rightarrow
T^*Q$ constructed from $L$ (no Riemannian structure is in
principle available on $Q$). In the case $\FL$ is a diffeomorphism,
the Lagrangian is called hyper-regular, and one can pull-back the
canonical symplectic form from $T^*Q$ to $TQ$ and define a
Hamiltonian system on $TQ$ (see \cite{Lewisblock} for details).
Given an element $\xi\in\g$, the locked Lagrangian $L_\xi\in
C^\infty (Q)$ is defined as $$ L_\xi (x)=L(\xi_Q(x)).$$ Also, the
locked momentum map is defined as the map $I_\xi:Q\rightarrow
\g^*$ that satisfies $$\langle
I_\xi(x),\eta\rangle=\frac{d}{dt}\restr{t=0}L_{\xi+t\eta}(x)$$ for
every $\eta\in\g$. We will also need the definition of the space
of admissible configuration variations at a point $x\in Q$, which
is \be\label{admissiblevariations}\mathcal{D}=\{\delta q\in T_xQ\,
:\,\D I_\xi\cdot\delta q\in\g_x^\circ\} ,\ee for a fixed element
$\xi\in\g$. Finally, given $x\in Q$ and $\xi\in\g$ the linearized
momentum map $I_x:\g\rightarrow \g^*$ is defined as
$$I_x(\eta)=\frac{d}{dt}\restr{t=0}I_{\xi+t\eta}(x)$$ for every
$\xi,\,\eta\in\g$. We will denote  its
generalized inverse by
$\tilde{I}_x^{-1}:\mathrm{range}\,I_x\rightarrow \g/\ker I_x$.

The Lagrangian Block-Diagonalization method gives sufficient
conditions for formal stability of relative equilibria in
Lagrangian systems. Under some assumptions, formal stability
implies $G_\mu$-stability. Here we shall not be concerned with
those differences, in order to provide a comparison between the
Lagrangian Block-Diagonalization and our singular reduced
Energy-Momentum Method. Furthermore, one needs a technical
condition relating $I_\xi$ and $I_x$ at a relative equilibrium in
order to be able to apply the method (see (3.15) in
\cite{Lewisblock}). However, for the particular case of Lagrangian
systems defining simple mechanical systems this is automatically
satisfied and thus it is not necessary for our comparison
objective. So we assume that the Lagrangian Block Diagonalization
method is applicable in order to simplify the exposition.

The next proposition collects the results of \cite{Lewisblock}
concerning relative equilibria and their $G_\mu$-stability. Let
$\mathbf{g}$ be the symmetric bilinear form on $T_xQ$ induced by
the hyper-regular Lagrangian $L$ and defined by $$\mathbf{g}
(v_x,w_x)=\frac{d}{dt}\restr{t=0}\frac{d}{ds}\restr{s=0}L(tv_x+sw_x).$$
By the hyper-regularity hypothesis of $L$, $\mathbf{g}(\cdot
,\cdot)$ is non-degenerate.
\begin{prop}[Lagrangian Block-Diagonalization, \cite{Lewisblock}]\label{lagstability}
Let $L\in C^G(TQ)$ be a hyper-regular Lagrangian invariant under
the tangent lifted action of the Lie group $G$ on $Q$. A point
$v_x\in TQ$ is a relative equilibrium for the Lagrangian system
defined on $TQ$ by $L$ if there is an element $\xi\in\g$ such that
$v_x=\xi_Q(x)$ and $x$ is a critical point of $L_\xi$. If the
bilinear form $$ \mathcal{B}=\left(-\ed_x^2L_\xi+\langle \D I_\xi
\cdot(\cdot),\tilde{I_x}^{-1}\left(\D I_\xi
\cdot(\cdot)\right)\rangle\right)\restr{\mathcal{D}}$$ is positive
(negative) semi-definite with kernel $\g_\mu\cdot x$ and
$\mathbf{g}_{\vert (\g\cdot x)^\perp}$ is positive (negative)
definite then the relative equilibrium is $G_\mu$-stable.
\end{prop} We now study Lagrangians defining simple mechanical
systems. It is a standard fact that if we are given a
$G$-invariant Riemannian metric $\ll\cdot,\cdot\gg$ on $Q$ and a
$G$-invariant function $V\in C^G(Q)$, then the Lagrangian
formulation of the associated simple mechanical system
\eqref{smshamiltonian} is \be\label{lagsms} L(v_x)=\frac 12
\parallel v_x\parallel^2- V(x).\ee
For $L$ of the form \eqref{lagsms} it is straightforward to
compute
$$\begin{array}{ll}
L_\xi & = -V_\xi\\
I_\xi (x) & = \II (x)(\xi)\\
I_x & = \II (x)\\
\mathbf{g} & =\ll\cdot,\cdot\gg .
\end{array}$$
Therefore, relative equilibria of the simple mechanical system
\eqref{lagsms} are defined by a velocity $\xi\in\g$ and a critical
point $x\in Q$ of $V_\xi$, recovering the well-know result for
critical points of the augmented Hamiltonian in simple mechanical
systems.

To see that the stability conditions of
Proposition~\ref{lagstability} are then equivalent to those given
in Theorem~\ref{REM}, one only needs to prove that at a relative
equilibrium, the space $\Sigma\subset T_xQ$ of
Definition~\ref{correction} is indeed a complement to $\g_\mu\cdot
x$ in $\mathcal{D}$, but this follows from Lemma~\ref{admissible}
and it has been already used in the proof of
Proposition~\ref{rigidinternalisos}. This shows that our stability
result, Theorem~\ref{REM}, is a Hamiltonian version of the
Lagrangian Block-Diagonalization method applied to simple
mechanical systems. In the same way, one can see that the extra
block-diagonalization construction carried out in subsection 3.3
of \cite{Lewisblock} is a consequence of the splitting of $V_s$ in
rigid and internal subspaces of
Section~\ref{singulararnoldsection} in this paper.
\paragraph{\bf{The normal form for the symplectic matrix.}}$\quad$\\
Besides the convenient form for
${\ed_{p_x}^2h_{\xi^\perp}}\rrestr{V_s}$ given in
Proposition~\ref{blockdiagonalforms} which allowed us to refine
Theorem~\ref{REM} and re-express it as
Corollary~\ref{blockdiagonalREM}, it is important to note the
particular expression for the symplectic matrix $\Omega$ of the
symplectic normal space $V_s$ identified with
$\q^\mu\oplus\Sigma_\text{int}\oplus\Sl^*$. In the free case (see
\cite{MarLec} and \cite{SiLeMar}), the explicit form for the
symplectic matrix, together with the one for
${\ed_{p_x}^2h_{\xi^\perp}}\rrestr{V_s}$, allows the authors to
obtain the linearization of the Hamiltonian vector field at a
relative equilibrium. This is an important observation, since the
study of this linearized vector field has applications in the
spectral and linear stability (and instability) of the relative
equilibrium under study, as well as for the identification of
possible bifurcations from it.

It seems that despite the generalization to the Lagrangian side
and non-free actions in \cite{Lewisblock} of the stability results
provided by the reduced Energy-Momentum Method of \cite{SiLeMar},
this feature has not been studied in detail for group actions with
singularities and the expression for the symplectic matrix
obtained in Proposition~\ref{blockdiagonalforms} for relative
equilibria with non-discrete configuration isotropies cannot be
found in the literature. Here we prove that our expression
\eqref{blocksigma} coincides in the regular (free) case with the
one obtained in \cite{SiLeMar}, equations (2.83) and (2.85).

Indeed (and if we assume that the Arnold form is non-degenerate)
given two e\-le\-ments
$(\eta_i,(\lambda^{a_i}_Q(x)+a_i),\beta_i)\in\q^\mu\oplus
w_\mathrm{int}\oplus \Sl^*\simeq V_s$, $i=1,2$, we will write
$\delta q_i=(\lambda^{a_i})_Q(x)+a_i$ and $\alpha (\delta
q_i)=\rII^{-1}(\J (\FL (\delta q_i)))=\lambda^{a_i}$. Then we have
the following:
\begin{prop}\label{finalmatrix} The expression for the symplectic matrix
$\Omega$ of Proposition~\ref{blockdiagonalforms} is equivalent to
$$\begin{array}{lll} \Omega
((\eta_1,\delta q_1,\beta_1),(\eta_2,\delta q_2,\beta_2)) & = &
\langle \mu, -[\eta_1,\eta_2]-[\eta_1,\alpha (\delta
q_2)]+[\eta_2,\alpha (\delta q_1)]\rangle\vspace{1mm}\\ & &
+\langle \beta_2,\delta q_1\rangle-\langle \beta_1,\delta
q_2\rangle -\ed \chi^{\xi^\perp}(x)(\delta q_1,\delta q_2),
\end{array}$$
with $\chi^{\xi^\perp}$ defined in \eqref{chi}. This coincides in
the regular case with equation (2.83) in \cite{SiLeMar}.
\end{prop}
\begin{proof}
Recall from Proposition~\ref{blockdiagonalforms} that we have
$$\begin{array}{lll} \Omega
((\eta_1,\delta q_1,\beta_1),(\eta_2,\delta q_2,\beta_2))& = &
\langle\mu, -[\eta_1,\eta_2]
-[\eta_1,\lambda^{a_2}]+[\eta_2,\lambda^{a_1}]\rangle\\ & &  +
\langle \beta_2,a_1\rangle-\langle
\beta_1,a_2\rangle-\langle\mu,[\lambda^{a_1},\lambda^{a_2}]\rangle\\
& &
 + \ll C(a_2)(\xi^\perp),a_1\gg_\Sl-\ll
C(a_1)(\xi^\perp),a_2\gg_\Sl .\end{array}$$
 Now note from the
definition of $\alpha$, that $[\eta_i,\alpha (\delta
q_j)]=[\eta_i,\lambda^{a_j}]$. Also, since
$\beta_i\in\Sl^*=(\g\cdot x)^\circ,\, i=1,2$, one has that
$\langle \beta_i,\delta q_j\rangle=\langle\beta_i,a_j\rangle$.
Therefore, the proposition will be proved if we show that
$$\ed
\chi^{\xi^\perp}(x)(\delta q_1,\delta
q_2)=\langle\mu,[\lambda^{a_1},\lambda^{a_2}]\rangle
 - \ll C(a_2)(\xi^\perp),a_1\gg_\Sl+\ll
C(a_1)(\xi^\perp),a_2\gg_\Sl.$$ To see this, we choose local
extensions $X=(\lambda^{a_1})_Q+\overline{a}_1$ and
$Y=(\lambda^{a_2})_Q+\overline{a}_2$ of $\delta q_1$ and $\delta
q_2$ near $x$ and then we use the formula for the exterior
derivative
$$\ed \chi^{\xi^\perp}(X,Y)=X
\left(\chi^{\xi^\perp}(Y)\right)-Y\left(\chi^{\xi^\perp}(X)\right)-\chi^{\xi^\perp}([X,Y]).$$
 We will need the following lemma, proved in \cite{PerRoSD2}, which shows
additional properties of the family of local vector fields
introduced in Section 3.
\begin{lemma}\label{normalformlemma}
The local vector fields defined by \eqref{localv}  satisfy
\begin{itemize}
\item[(i)] $\left[\overline{v}_a,\overline{v}_b\right](x)=0$ for
every $v_a,\,v_b\in\Sl$ \item[(ii)]
$\left[\lambda_Q,\overline{v}\right](x)=
  0$  for every  $\lambda\in\lir,\,v\in\Sl$
 \item[(iii)] There is a small enough open
neighbourhood $O\ni x$ such that for every $x'\in G\cdot x\cap
O,\,v\in\Sl$ and $\xi\in\g$, $\ll\overline{v}(x'),\xi_Q(x')\gg
=0.$ \item[(iv)] Let $x'=\sigma([g,s])$ with $\sigma$ defined in
\eqref{tube} and $\lambda\in\lir$, then  $$\lambda_Q(x')=\frac
{d}{dt}\restr{t=0}\sigma([(\exp_e t\lambda)g,s]).$$
\end{itemize}
\end{lemma}
Recall the definition $\chi^{\xi^\perp}=\FL (\xi^\perp_Q)$ and let
us compute separately the group orbit and slice contributions for
$\ed \chi^{\xi^\perp}(x)(X,Y)$ using the usual properties of the
Levi-Civita connection.
$$\begin{array}{ll}
\ed \chi^{\xi^\perp}(x)(a_1,a_2) & =
\overline{a}_1\left(\ll\xi^\perp_Q,\overline{a}_2\gg\right)(x)-
\overline{a}_2\left(\ll\xi^\perp_Q,\overline{a}_1\gg\right)(x)\vspace{1mm}\\
& =
\ll\nabla_{\overline{a}_1}\xi^\perp_Q(x),a_2\gg+\ll\xi^\perp_Q(x),\nabla_{\overline{a}_1}\overline{a}_2(x)\gg\vspace{1mm}\\
& -\ll\nabla_{\overline{a}_2}\xi^\perp_Q(x),a_1\gg-
\ll\xi^\perp_Q(x),\nabla_{\overline{a}_2}\overline{a}_1(x)\gg\vspace{1mm}\\
& =\ll C(a_1)(\xi^\perp),a_2\gg_\Sl-\ll
C(a_2)(\xi^\perp),a_1\gg_\Sl ,
\end{array}$$
where we have used the definition of $C$ \eqref{defC}, together
with item (v) in Theorem \ref{covariant} and the fact that
$$\begin{array}{l}\ll\xi^\perp_Q(x),\nabla_{\overline{a}_1}\overline{a}_2(x)\gg-
\ll\xi^\perp_Q(x),\nabla_{\overline{a}_2}\overline{a}_1(x)\gg\vspace{1mm}
\\ \hspace{1.5cm} =\ll\xi^\perp_Q(x),T(\overline{a}_1,\overline{a}_2)(x)\gg +
\ll\xi^\perp_Q(x),[\overline{a}_1,\overline{a}_2](x)]\gg=0\end{array}$$
by item (i) of Lemma~\ref{normalformlemma} as well as noting
that the Levi-Civita connection has zero torsion
($T(X,Y)=\nabla_XY-\nabla_YX-[X,Y]=0$, for every pair of vector
fields $X,Y$). In the same way we can compute
$$\begin{array}{ll}
\ed \chi^{\xi^\perp}(x)((\lambda^{a_i})_Q(x),a_j) & =
(\lambda^{a_i})_Q\left(\ll\xi^\perp_Q,\overline{a}_j\gg\right)(x)
-\overline{a}_j\left(\ll\xi^\perp_Q,(\lambda^{a_i})_Q\gg\right)(x)\vspace{1mm}\\
&\hspace{5mm} -
\chi^{\xi^\perp}\left(\left[(\lambda^{a_i})_Q,\overline{a}_j\right]\right)(x)
= - (\D\II\cdot a_j)(\xi^\perp,\lambda^{a_i}),
\end{array}$$
since by (iii) and (ii) in Lemma~\ref{normalformlemma} the first
and last contributions vanish. Finally
$$\begin{array}{l}
\ed \chi^{\xi^\perp}(x)((\lambda^{a_1})_Q(x),(\lambda^{a_2})_Q(x))
=\vspace{1mm}\\ \hspace{1cm}
(\lambda^{a_1})_Q\left(\ll\xi^\perp_Q,(\lambda^{a_2})_Q\gg\right)(x)-
(\lambda^{a_2})_Q\left(\ll\xi^\perp_Q,(\lambda^{a_1})_Q\gg\right)(x)\vspace{1mm}\\
\hspace{1cm}
+\ll\xi^\perp_Q(x),([\lambda^{a_1},\lambda^{a_2}])_Q(x)\gg =
(\D\II\cdot(\lambda^{a_1})_Q(x))(\xi^\perp,\lambda^{a_2})\vspace{1mm}\\
\hspace{1cm} -
(\D\II\cdot(\lambda^{a_2})_Q(x))(\xi^\perp,\lambda^{a_1})+\langle\mu,[\lambda^{a_1},\lambda^{a_2}]\rangle
.
\end{array}$$
Putting together all the contributions we obtain
$$\begin{array}{ll}
\ed \chi^{\xi^\perp}(x)(\delta q_1,\delta q_2) &
=\langle\mu,[\lambda^{a_1},\lambda^{a_2}]\rangle
+\ll C(a_1)(\xi^\perp),a_2\gg_\Sl-\ll C(a_2)(\xi^\perp),a_1\gg_\Sl\vspace{1mm}\\
& +(\D\II\cdot\delta
q_1)(\xi^\perp,\lambda^{a_2})-(\D\II\cdot\delta
q_2)(\xi^\perp,\lambda^{a_1}) .
\end{array}$$
The last two terms of the above expression vanish since for
$i=1,2$, $\lambda^{a_i}\in\q^\mu\subset \mathfrak{t}$ and
$(\lambda^{a_i},a_i)\in w_\mathrm{int}$, and by its definition
\eqref{wintchar} $\Proj_{\mathfrak{t}}\,\left[(\D\II\cdot\delta
q_1)(\xi^\perp)\right]=0$. This finishes the proof.
\end{proof}
\paragraph{\bf{Acknowledgements.}} I would like to thank Esmeralda Sousa-Dias,
Matt Perlmutter and Tanya Schmah for many illuminating discussions
over the past years which have helped me during the realization of this
project. I would also like to thank CDS of the California Institute of Technology
for their hospitality and facilities  during a
research visit in November 2004 when part of this project was
carried out and in particular Jerrold Marsden for several
useful suggestions and comments on this research during that
period. This research was partially supported by the EU through
funds for the European Research Training Network ``MASIE'', vide 
contract HPRN-CT-2000-00113. Finally, I would like to thank the
referees for their comments on the manuscript.


\begin{thebibliography}{99}
\bibitem{Arnold2} V.I Arnold [1966], Sur la g\'eometrie differentielle
des groupes de Lie de dimension infinie et ses applications \`{a}
l'hydrodynamique des fluids parfaits, {\it Ann. Inst. Fourier,
Grenoble}, {\bf 16}, 319--361.



\bibitem{DuiKol} J.J. Duistermaat and J.A.C. Kolk [2000], {\it Lie
groups}, Universitext, Springer-Verlag.

\bibitem{Golu} M. Golubitsky, I. Stewart and D.G. Schaeffer [1988], Singularities
and groups in bifurcation theory. Vol. II. Applied Mathematical
Sciences, 69. Springer-Verlag, New York.

\bibitem{Kara} A. V. Karapetyan [1988], The Routh Theorem and its
extensions, {Colloq. Math. Soc. Janos Bolyai} {\bf 53}, {\it
Qualitative Theory of Differential Equations}, 271--290.

\bibitem{LeSi} E. Lerman and S.F. Singer [1998], Stability and
persistence of relative equilibria at singular values of the
moment map, {\it Nonlinearity} {\bf 11}, 1637--1649.

\bibitem{Lewisblock} D. Lewis [1992], Lagrangian
Block-Diagonalization, {\it Journal of Dynamics and Differential
Equations}, vol. {\bf 4}, No. {\bf 1}, 1--41.

\bibitem{Lewisdrop} D. Lewis [1993], Bifurcations of liquid drops. {\it Nonlinearity} {\bf 6}, 491-522.


\bibitem{LeRaSiMar} D. Lewis, T.S. Ratiu, J.C. Simo and J.E.
Marsden [1992], The heavy top: a geometric treatment, {\it
Nonlinearity} {\bf 5}, 1--48.


\bibitem{MarLec} J.E. Marsden [1992], {\it Lectures on Mechanics}, Lecture Note Series
{\bf 174}, LMS, Cambridge University Press.

\bibitem{Montaldi} J. Montaldi [1997], Persistence and stability of
relative equilibria, {\it Nonlinearity} {\bf 10}, 449--466.


\bibitem{OrRa99} J.-P. Ortega  and T.S. Ratiu [1999], Stability of
Hamiltonian relative equilibria, {\it Nonlinearity} {\bf 12},
693--720.

\bibitem{Palais} R.S. Palais [1961], On the existence of slices
for actions of non-compact Lie groups, {\it Ann. of Math.}, {\bf
73}, 265--323.

\bibitem{Paternain} G.P. Paternain [1999], {\it Geodesic Flows}, Progress in Mathematics {\bf 180}, Birkhauser.

\bibitem{Pascal} M. Pascal [1975], Sur la recherche des mouvements
stationnaires dans les systemes ayant des variables cycliques,
{\it Celestial Mechanics}, {\bf 12}, 337--358.



\bibitem{Patrick} G.W. Patrick [1992], Relative equilibria in Hamiltonian systems: the dynamic
interpretation of nonlinear stability on a reduced phase space,
{\it J. Geom. Phys.} {\bf 9}, no. 2, 111--119.

\bibitem{PatRobWul} G.W. Patrick, M. Roberts and C. Wulff [2002],
Stability of Poisson Equilibria and Hamiltonian Relative
Equilibria by Energy Methods, math.DS/0201239.


\bibitem{PerRoSD2} M. Perlmutter, M Rodr\'{\i}guez-Olmos and E.
Sousa-Dias [2004], The symplectic normal space of a
cotangent-lifted action, {\it Diff. Geom. App.} to appear. Preprint math.SG/0501207.

\bibitem{Routh} E.J. Routh [1877], {\it A Treatise on the
Stability of a Given State of Motion}, MacMillan and Co, London.

\bibitem{SiLeMar} J.C. Simo, D. Lewis and J.E. Marsden [1991], Stability
of relative equilibria. Part I: The reduced Energy-Momentum
Method, {\it Arch. Rational Mech. Anal.}, {\bf 115}, 15--59.

\bibitem{Smale} S. Smale [1970], Topology and Mechanics I, {\it
Inventiones Math.}, {\bf 10}, 305--331.


\end{thebibliography}
\end{document}